\documentclass[review,onefignum,onetabnum]{siamart171218}
\usepackage{lipsum}
\usepackage{amsfonts}
\usepackage{graphicx}
\usepackage{epstopdf}
\usepackage{algorithmic}

\usepackage{url}

\usepackage{amssymb,amsmath}
\usepackage{amsfonts}
\usepackage{comment}
\usepackage{color}
\usepackage{mathrsfs}

\ifpdf
  \DeclareGraphicsExtensions{.eps,.pdf,.png,.jpg}
\else
  \DeclareGraphicsExtensions{.eps}
\fi

\newsiamremark{remark}{Remark}
\newsiamremark{hypothesis}{Hypothesis}
\crefname{hypothesis}{Hypothesis}{Hypotheses}
\newsiamthm{claim}{Claim}

\headers{Yaglom limit for Stochastic Fluid Models}{N.G. Bean, M. O'Reilly, and Z. Palmowski}

\title{Yaglom limit for Stochastic Fluid Models\thanks{
\funding{We would like to thank the Australian Research Council for funding this research through Linkage Project LP140100152.
Zbigniew Palmowski was partially supported by the National Science
Centre (Poland) under the grant 2018/29/B/ST1/00756.
}}}

\author{
Nigel G. Bean\thanks{Australian Research Centre of Excellence for Mathematical and Statistical Frontiers.
School of Mathematical Sciences, University of Adelaide, SA 5005, Australia(\email{nigel.bean@adelaide.edu.au}).}
\and 		
Ma{\l}gorzata M. O'Reilly
\thanks{Australian Research Centre of Excellence for Mathematical and Statistical Frontiers. Discipline of Mathematics, University of Tasmania, Hobart TAS 7001, Australia (\email{malgorzata.oreilly@utas.edu.au}).}
\and
Zbigniew Palmowski\thanks{Faculty of Pure and Applied Mathematics, Wroc{\l}aw University of Science and Technology, ul. Wybrze\.{z}e Wyspia\'{n}skiego 27,
50-370 Wroc{\l}aw, Poland (\email{zbigniew.palmowski@pwr.edu.pl}).}
}

\usepackage{amsopn}

\newtheorem{cor}{Corollary}
\numberwithin{cor}{section}

\newcommand{\bpsi}{\mbox{\boldmath$\psi$}}
\newcommand{\bxi}{\mbox{\boldmath$\xi$}}

\newcommand{\bphi}{\mbox{\boldmath$\phi$}}
\newcommand{\bmu}{\mbox{\boldmath$\mu$}}

\newtheorem{example}{Example}

\ifpdf
\hypersetup{
  pdftitle={Yaglom limit for the Stochastic Fluid Models},
  pdfauthor={N.G. Bean, M. O'Reilly, and Z. Palmowski}
}
\fi

\begin{document}

\maketitle

\begin{abstract}
  In this paper we provide the analysis of the limiting conditional distribution (Yaglom limit) for stochastic fluid models (SFMs), a key class of models in the theory of matrix-analytic methods.

So far, transient and stationary analyses of the SFMs have been only considered in the literature. The limiting conditional distribution gives useful insights into what happens when the process has been evolving for a long time, given its busy period has not ended yet.

We derive expressions for the Yaglom limit in terms of the singularity~$s^*$ such that the key matrix of the SFM, ${\bf\Psi}(s)$, is finite (exists) for all $s\geq s^*$ and infinite for $s<s^*$. We show the uniqueness of the Yaglom limit and illustrate the application of the theory with simple examples.
\end{abstract}



\section{Introduction}

Let $\{ (\varphi(t)):t\geq 0\}$ be an irreducible, positive-recurrent, continuous-time Markov chain (CTMC) with some finite state space $\mathcal{S}=\{1,2,\ldots ,n\}$ and infinitesimal generator ${\bf T}$. Let $\{(\varphi(t),X(t)):t\geq 0\}$ be a Markovian stochastic fluid model (SFM)~\cite{AR,AR2,AR3,Asm,BOT2,BOT,BOT5,rama96,rama99}, with phase variable $\varphi(t)\in\mathcal{S}$, level variable $X(t)\geq 0$, and constant rates $c_i\in\mathbb{R}$, for all $i\in\mathcal{S}$. The model assumes that when $\varphi(t)=i$ and $X(t)>0$, then the rate at which the level is changing is $c_i$, and when $\varphi(t)=i$ and $X(t)=0$, then the rate at which the level is changing is $max\{0,c_i\}$. Therefore, we refer to the CTMC $\{ (\varphi(t)):t\geq 0\}$ as the process that is driving (or modulating) the SFM $\{(\varphi(t),X(t)):t\geq 0\}$.

SFMs are a key class of models in the theory of matrix-analytic methods~\cite{Lat-Ram,LR2013,Qi2013}, which comprises methodologies for the analysis of Markov chains and Markovian-modulated models, that lead to efficient algorithms for numerical computation.

Let $\mathcal{S}_1=\{i\in\mathcal{S}:c_i>0\}$, $\mathcal{S}_2=\{i\in\mathcal{S}:c_i<0\}$, $\mathcal{S}_0=\{i\in\mathcal{S}:c_i=0\}$, and partition the generator as
\begin{displaymath}
{\bf T}=\left[
\begin{array}{ccc}
{\bf T}_{11} & {\bf T}_{12} &  {\bf T}_{10}\\
{\bf T}_{21} & {\bf T}_{22} & {\bf T}_{20}\\
{\bf T}_{01} & {\bf T}_{02} & {\bf T}_{00}
\end{array}
\right],
\end{displaymath}
according to $\mathcal{S}=\mathcal{S}_1\cup\mathcal{S}_2\cup\mathcal{S}_0$.

We assume that the process is stable, that is
\begin{eqnarray}
\mu = \sum_{i\in\mathcal{S}} c_i\xi_i &<& 0,
\end{eqnarray}
where $\bxi=[\xi_i]_{i\in\mathcal{S}}$ is the stationary distribution vector of the Markov chain $\{(\varphi(t)):t\geq 0\}$.

So far, the analysis of SFMs has focused on the transient and stationary behaviour. In this paper, we are interested in the behaviour of the process conditional on absorption not having taken place; where absorption means that the busy period of the process has ended, that is, the process has not hit the level zero as yet. For $x\geq 0$, let $\theta(x)=\inf\{t>0:X(t)=x\}$ be the first time at which the process reaches level $x$. To this end, we define the following quantity, referred to as the Yaglom limit.
\begin{definition}\label{def:mu}
Define the matrix $\bmu(dy)^{(x)}=[\mu(dy)^{(x)}_{ij}]_{i,j\in\mathcal{S}}$, $x,y> 0$, such that,
\begin{equation}
\mu(dy)^{(x)}_{ij}=\lim_{t\to\infty}P(X(t)\in dy,\varphi(t)=j\ |\ \theta(0)>t,X(0)=x,\varphi(0)=i),
\end{equation}
and matrix $\bmu(dy)^{(0)}=[\mu(dy)^{(0)}_{ij}]_{i\in\mathcal{S}_1,j\in\mathcal{S}}$, $y> 0$, such that
\begin{equation}
\mu(dy)^{(0)}_{ij}=\lim_{t\to\infty}P(X(t)\in dy,\varphi(t)=j\ |\ \theta(0)>t,X(0)=0,\varphi(0)=i),
\end{equation}
whenever the limit exists.
We refer to $\mu(dy)^{(x)}_{ij}$ as the limiting conditional distribution (Yaglom limit) of observing the process in level $y$ and phase $j$, given the process started from level $x$ in phase $i$ at time zero, and has been evolving without hitting level zero.
\end{definition}
\begin{remark}
In general, for Markov processes there are no sufficient conditions that we can refer to under which there exists
Yaglom limit or quasi-stationary distribution .
Usually, the existence of Yaglom limit is proved case by case. Here, we prove that it exists for our model.
\end{remark}

We partition $\bmu(dy)^{(x)}$, $x>0$, according to $\mathcal{S}_1\cup\mathcal{S}_2\cup\mathcal{S}_0 \times \mathcal{S}_1\cup\mathcal{S}_2\cup\mathcal{S}_0$ as
\begin{equation}
\bmu(dy)^{(x)}=
\left[
\begin{array}{ccc}
\bmu(dy)^{(x)}_{11}&\bmu(dy)^{(x)}_{12}&\bmu(dy)^{(x)}_{10}\\
\bmu(dy)^{(x)}_{21}&\bmu(dy)^{(x)}_{22}&\bmu(dy)^{(x)}_{20}\\
\bmu(dy)^{(x)}_{01}&\bmu(dy)^{(x)}_{02}&\bmu(dy)^{(x)}_{00}\\
\end{array}
\right],
\end{equation}
and partition its row sums accordingly, as
\begin{equation}
\bmu(dy)^{(x)}{\bf 1}=
\left[
\begin{array}{c}
\bmu(dy)^{(x)}_{1}\\
\bmu(dy)^{(x)}_{2}\\
\bmu(dy)^{(x)}_{0}\\
\end{array}
\right],
\end{equation}
where ${\bf 1}$ denotes a vector of ones of appropriate size, so that $\bmu(dy)^{(x)}_{1}=\bmu(dy)^{(x)}_{11}{\bf 1}+\bmu(dy)^{(x)}_{12}{\bf 1}+\bmu(dy)^{(x)}_{10}{\bf 1}$, and so on.

We partition $\bmu(dy)^{(0)}$ according to $\mathcal{S}_1\times\mathcal{S}_1\cup\mathcal{S}_2\cup\mathcal{S}_0$ as
\begin{equation}
\bmu(dy)^{(0)}=
\left[
\begin{array}{ccc}
\bmu(dy)^{(0)}_{11}&\bmu(dy)^{(0)}_{12}&\bmu(dy)^{(0)}_{10}\\
\end{array}
\right],
\end{equation}
and let
\begin{equation}
\bmu(dy)^{(0)}_{1}=\bmu(dy)^{(0)}{\bf 1}.
\end{equation}

This paper is the first analysis of the Yaglom limit of SFMs. We derive expressions for the Yaglom limit, show its uniqueness and illustrate the theory with simple examples.
Yaglom limit concerns some Markov process $X(t)$ and some finite a.s. absorption time $\theta$
(usually, first exit time from some set such as a positive half-line), and is defined by
\begin{equation}
\mu(dy)=\lim_{t\rightarrow +\infty} P(X(t)\in dy\ |\ \theta >t).
\end{equation}
It describes the state of the Markov system conditioned on surviving killing coming from~$\theta$ for a very long time.
Yaglom limit is strongly related with so-called quasi-stationary distribution that satisfies
\begin{equation}
P_\mu(X(t)\in dy|\theta >t)=\mu(dy),
\end{equation}
see for example \cite{DS1965}. In particular, the Yaglom limit $\mu$ (if exists) is necessarily quasi-stationary
but it may be difficult to show its uniqueness~\cite[Section 3]{bookAsm}. In other words, there might be more quasi-stationary laws and Yaglom limit might be one of them.
It might be the case as well that there exists quasi-stationary distribution but Yaglom is not well-defined.

A related class of models in the theory of matrix-analytic methods, is Quasi-Birth-and-Death process (QDBs)~\cite{Lat-Ram}, in which the level variable is discrete.
The quasi-stationary analysis of the QBDs has been provided in~\cite{BBL1997,BPT1998,BPT2000}, along with several examples of areas of applications, which are relevant here as well, due to the similar application potential of the QBDs and SFMs~\cite{BO2012uni}.

Information on quasi-stationary distributions (QS) for other
Markov processes can be found in the classical works of Seneta and Vere-Jones \cite{MR0207047}, Tweedie \cite{MR0443090}, Jacka and Roberts \cite{MR1363332}.
The bibliographic database of Pollet \cite{pollett} gives detailed history of
quasi-stationary distributions. In particular,
Yaglom \cite{MR0022045} was the first to explicitly  identify QS
distributions for the subcritical Bienaym\'e-Galton-Watson branching process.
Part of the results on QS distributions
concern Markov chains on positive integers with an absorbing state at the origin
\cite{MR2986807, MR1334159, MR0346932, MR0207047, MR1133722,MR3247530}.
Other objects of study are the extinction probabilities for continuous-time branching process and the Fleming-Viot process \cite{MR3498004, MR2318407, MR2299923}.
A separate topic is the L\'evy processes exiting from the positive half-line or a cone.
Here the case of the Brownian motion with drift was resolved by Martinez and San Martin \cite{MR1303922}, complementing the result for random walks obtained by Iglehart \cite{MR0368168}.
The case of more general L\'evy processes was studied by \cite{jacone, MR2248228, MR0292201, MR2959448}.
One-dimensional self-similar processes, including the  symmetric $\alpha$-stable  L\'evy process,  were subject of interest of \cite{MR2971725}.

The rest of the paper is structured as follows. In Section~\ref{sec:LSTs} we define the Laplace-Stieltjes Transforms (LSTs) which form the key building blocks of the analysis and in Section~\ref{sec:heaviside} we outline the approach based on the Heaviside principle. The key results of this paper are contained in Section~\ref{sec:applH}. To illustrate the theory we construct a simple example with scalar parameters, which we analyse throughout the paper, as we introduce the theory. In Section~\ref{sec:s_ex} we analyse another example, with matrix parameters, where we provide some numerical output as well.

\section{The Laplace-Stieltjes Transforms}
\label{sec:LSTs}

Note that by Definition~\ref{def:mu}, for $x\geq 0$,
\begin{eqnarray}
\mu(dy)^{(x)}_{ij}
&=&
\lim_{t\to\infty}P(X(t)\in dy,\varphi(t)=j\ |\ \theta(0)>t,X(0)=x,\varphi(0)=i)\nonumber\\
&=&
\lim_{t\to\infty}
\frac{
P(X(t)\in dy,\varphi(t)=j,\theta(0)>t\ |\
X(0)=x,\varphi(0)=i)
}
{
P(\theta(0)>t\ |\ X(0)=x,\varphi(0)=i)
},
\label{mu_fraction}
\end{eqnarray}
and, for all $x\geq 0$ and $y>0$, define the matrix ${\bf E}(dy)^{(x)}(s)=[E(dy)^{(x)}_{ij}(s)]_{i,j\in\mathcal{S}}$ and the vector ${\bf E}^{(x)}(s)=[E^{(x)}_i(s)]_{i\in\mathcal{S}}$ , which record the corresponding Laplace-Stieltjes Transforms (LSTs),
\begin{eqnarray}
E(dy)^{(x)}_{ij}(s)&=&\int_0^\infty E(e^{-s t} {\bf 1}
\{X(t)\in dy,\varphi(t)=j,\theta(0)>t\}\ |\
X(0)=x,\varphi(0)=i) dt,\nonumber\\
E^{(x)}_i(s)&=&\int_0^\infty E(e^{-s t} {\bf 1}\{\theta(0)>t\}\ |\ X(0)=x,\varphi(0)=i)dt,
\end{eqnarray}
where ${\bf 1}\{\cdot\}$ denotes an indicator function. We have,
\begin{eqnarray}\label{Edy}
{\bf E}^{(x)}(s)&=&\int_{y=0}^{\infty}{\bf E}(dy)^{(x)}(s){\bf 1}
.
\end{eqnarray}

We  partition ${\bf E}(dy)^{(x)}(s)$, $x>0$, according to $\mathcal{S}\times \mathcal{S}$ for $\mathcal{S}=\mathcal{S}_1\cup\mathcal{S}_2\cup\mathcal{S}_0$ as
\begin{equation}
{\bf E}(dy)^{(x)}(s)=
\left[
\begin{array}{ccc}
{\bf E}(dy)^{(x)}(s)_{11}&{\bf E}(dy)^{(x)}(s)_{12}&{\bf E}(dy)^{(x)}(s)_{10}\\
{\bf E}(dy)^{(x)}(s)_{21}&{\bf E}(dy)^{(x)}(s)_{22}&{\bf E}(dy)^{(x)}(s)_{20}\\
{\bf E}(dy)^{(x)}(s)_{01}&{\bf E}(dy)^{(x)}(s)_{02}&{\bf E}(dy)^{(x)}(s)_{00}\\
\end{array}
\right],
\end{equation}
and ${\bf E}^{(x)}(s)$, $x>0$, as
\begin{equation}
{\bf E}^{(x)}(s)=
\left[
\begin{array}{c}
{\bf E}^{(x)}(s)_1\\
{\bf E}^{(x)}(s)_2\\
{\bf E}^{(x)}(s)_0
\end{array}
\right].
\end{equation}

We  partition ${\bf E}(dy)^{(0)}(s)$ according to $\mathcal{S}_1\times\mathcal{S}_1\cup\mathcal{S}_2\cup\mathcal{S}_0$ as
\begin{equation}
{\bf E}(dy)^{(0)}(s)=
\left[
\begin{array}{ccc}
{\bf E}(dy)^{(0)}(s)_{11}&{\bf E}(dy)^{(0)}(s)_{12}&{\bf E}(dy)^{(0)}(s)_{10}\\
\end{array}
\right]
\end{equation}
and let
\begin{equation}
{\bf E}(dy)^{(0)}(s)_{1}={\bf E}(dy)^{(0)}(s){\bf 1}.
\end{equation}

Denote ${\bf C}_1=diag(c_i)_{i\in\mathcal{S}_1}$, ${\bf C}_2=diag(|c_i|)_{i\in\mathcal{S}_2}$, and let ${\bf Q}(s)$ be the key fluid generator matrix ${\bf Q}(s)$ introduced in~\cite{BOT},
\begin{equation}
{\bf Q}(s)=
\left[
\begin{array}{cc}
{\bf Q}_{11}(s) & {\bf Q}_{12}(s) \\
{\bf Q}_{21}(s) & {\bf Q}_{22}(s)
\end{array}
\right],
\end{equation}
where the block matrices are given by,
\begin{eqnarray}
{\bf Q}_{22}(s)&=&
{\bf C}^{-1}_{2}
\left(
{\bf T}_{22}-s{\bf I}
-
{\bf T}_{20}
({\bf T}_{00}-s{\bf I})^{-1}{\bf T}_{02}
\right),
\nonumber\\
{\bf Q}_{11}(s)&=&
{\bf C}^{-1}_{1}
\left(
{\bf T}_{11}-s{\bf I}
-
{\bf T}_{10}
({\bf T}_{00}-s{\bf I})^{-1}
{\bf T}_{01}\right),
\nonumber\\
{\bf Q}_{12}(s)&=&
{\bf C}^{-1}_{1}
\left(
{\bf T}_{12}
-
{\bf T}_{10}
({\bf T}_{00}-s{\bf I})^{-1}
{\bf T}_{02}\right)
,
\nonumber\\
{\bf Q}_{21}(s)&=&
{\bf C}^{-1}_{2}
\left(
{\bf T}_{21}
-
{\bf T}_{20}
({\bf T}_{00}-s{\bf I})^{-1}
{\bf T}_{01}\right)
\label{eq:Qs},
\end{eqnarray}
where ${\bf Q}(s)$ exists for all real $s$ such that $({\bf T}_{00}-s{\bf I})^{-1}=\int_{t=0}^{\infty}e^{-st}e^{{\bf T}_{00}t}dt<\infty$
or for all real $s$ when $S_0=\emptyset$.

Also, let ${\bf\Psi}(s)$ be the key matrix for SFMs~\cite{BOT} such that, for all $i\in\mathcal{S}_1, j\in\mathcal{S}_2$,
\begin{equation}
[{\bf\Psi}(s)]_{ij}=E(e^{-s\theta( 0)}
 {\bf 1}\{\theta(0)<\infty,\varphi(\theta(0))=j\}\ |\ \varphi(0)=i,X(0)=0)
\label{def_Psi_s}
\end{equation}
is the LST of the first return time to the original level $0$ and doing so in phase $j$, given start in level $0$ in phase $i$. Let $\bpsi(t)$ be the corresponding density so that ${\bf\Psi}(s)=\int_{t=0}^{\infty}e^{-st}\bpsi(t)dt$. Clearly, ${\bf\Psi}(s)>0$ for real $s$ such that ${\bf\Psi}(s)<\infty$ exists.

Define matrices
\begin{eqnarray}
	{\bf K}(s)&=&{\bf Q}_{11}(s)+{\bf\Psi}(s){\bf Q}_{21}(s),\nonumber\\
{\bf D}(s)&=&{\bf Q}_{22}(s)+{\bf Q}_{21}(s){\bf\Psi}(s),
\label{eq:KsDs}
\end{eqnarray}
and note that by the assumed stability of the process, the spectra of ${\bf K}(s)$ and $(-{\bf D}(s))$ are separate for $s\geq 0$ by~\cite{BOT2,BOT,BOT5}, that is, $sp({\bf K}(s))\cap sp(-{\bf D}(s))=\varnothing$.

We extend the result in~\cite[Equation (23)]{BOT} from $Re(s)\geq 0$ to all real $s$ such that ${\bf\Psi}(s)<\infty$ exists.
\begin{lemma}\label{elm:psi}
For all real $s$ such that ${\bf\Psi}(s)<\infty$ exists, the matrix ${\bf\Psi}(s)$ is a solution of the Riccati equation,
\begin{equation}\label{Ric}
{\bf Q}_{12}(s)+{\bf Q}_{11}(s){\bf X}+{\bf X}{\bf Q}_{22}(s)+{\bf X}{\bf Q}_{21}(s){\bf X}={\bf 0}.
\end{equation}
\end{lemma}
{\bf Proof:} Suppose $s$ is real and ${\bf\Psi}(s)<\infty$. Then, by~\cite[Theorem~1]{BOT} and~\cite[Algorithm~1 of Section 3.1]{BOT5},
\begin{equation}
\infty> {\bf\Psi}(s)=\int_{y=0}^{\infty}e^{{\bf Q}_{11}(s)y}
({\bf Q}_{12}(s)+{\bf\Psi}(s){\bf Q}_{21}(s){\bf\Psi}(s))
e^{{\bf Q}_{22}(s)y}dy.
\end{equation}
Then, letting
\begin{eqnarray}
{\bf\Psi}(s,y)&=& e^{{\bf Q}_{11}(s)y}
({\bf Q}_{12}(s)+{\bf\Psi}(s){\bf Q}_{21}(s){\bf\Psi}(s))
 e^{{\bf Q}_{22}(s)y},
\end{eqnarray}
gives
\begin{eqnarray}
\frac{\partial}{\partial y} {\bf\Psi}(s,y) &=&{\bf Q}_{11}(s){\bf\Psi}(s,y)+{\bf\Psi}(s,y){\bf Q}_{22}(s),\\
{\bf\Psi}(s,0)&=&{\bf Q}_{12}(s)+{\bf\Psi}(s){\bf Q}_{21}(s){\bf\Psi}(s),\\
\lim_{y\to\infty}{\bf\Psi}(s,y)&=&{\bf 0},
\end{eqnarray}
and so
\begin{eqnarray}
{\bf Q}_{11}(s){\bf\Psi}(s)+{\bf\Psi}(s){\bf Q}_{22}(s)&=&
\int_{y=0}^{\infty}
\frac{\partial}{\partial y}\left( {\bf\Psi}(s,y) \right)
dy\nonumber\\
&=&\lim_{y\to\infty}{\bf\Psi}(s,y)-{\bf\Psi}(s,0),
\end{eqnarray}
and
\begin{equation}
{\bf 0}={\bf Q}_{12}(s)+{\bf Q}_{11}(s){\bf\Psi}(s)+{\bf\Psi}(s){\bf Q}_{22}(s)+{\bf\Psi}(s){\bf Q}_{21}(s){\bf\Psi}(s),
\end{equation}
which implies that ${\bf\Psi}(s)$ is a solution of~\eqref{Ric}. \rule{9pt}{9pt}

Below, we state expressions for ${\bf E}(dy)^{(0)}(s)$ derived in~\cite[Theorem 3.1.1]{AR2} and~\cite[Theorem 2]{BOoperator}.
\begin{lemma}\label{lem:Edy_0}
We have
\begin{eqnarray}
{\bf E}(dy)^{(0)}(s)_{11}&=&e^{{\bf K}(s)y}
{\bf C}_1^{-1}dy,\label{eq:Edy011}\\
{\bf E}(dy)^{(0)}(s)_{12}&=&e^{{\bf K}(s)y}{\bf\Psi}(s)
{\bf C}_2^{-1}dy,\label{eq:Edy012}\\
{\bf E}(dy)^{(0)}(s)_{10}&=&
\left[
\begin{array}{cc}
{\bf E}(dy)^{(0)}(s)_{11}&{\bf E}(dy)^{(0)}(s)_{12}
\end{array}
\right]
\left[
\begin{array}{c}
{\bf T}_{10}\\
{\bf T}_{20}
\end{array}
\right]
\nonumber\\
&&
\times
(-({{\bf T}}_{00}-s{\bf I})^{-1})
dy.
\end{eqnarray}	
\end{lemma}

Next, we derive expressions for ${\bf E}(dy)^{(x)}(s)$, $x>0$.
The formal proof was already given by Ahn and Ramaswami \cite{SR_2006}.
Since it is a crucial lemma for whole further analysis we decided
to add its proof for completeness of all arguments.
\begin{lemma}\label{th:Edy}
For $x>0$ we have,
\begin{eqnarray}
{\bf E}(dy)^{(x)}(s)_{21}&=&
\int_{z=0}^{\min\{x,y\}}
e^{{\bf D}(s)(x-z)}{\bf Q}_{21}(s)e^{{\bf K}(s)(y-z)}{\bf C}_1^{-1}dzdy
,\nonumber\\
{\bf E}(dy)^{(x)}(s)_{22}&=&
{\bf E}(dy)^{(x)}(s)_{21}{\bf C}_1{\bf\Psi}(s){\bf C}_2^{-1}
+e^{{\bf D}(s)(x-y)}{\bf C}_2^{-1}{\bf 1}\{y<x\}dy,\nonumber\\
{\bf E}(dy)^{(x)}(s)_{20}&=&
\left[
\begin{array}{cc}
{\bf E}(dy)^{(x)}(s)_{21}&{\bf E}(dy)^{(x)}(s)_{22}
\end{array}
\right]
\left[
\begin{array}{c}
{\bf T}_{10}\\
{\bf T}_{20}
\end{array}
\right]
\nonumber\\
&&
\times
(-({{\bf T}}_{00}-s{\bf I})^{-1})
dy
,\nonumber
\\
{\bf E}(dy)^{(x)}(s)_{11}&=&{\bf\Psi}(s){\bf E}(dy)^{(x)}(s)_{21}
+e^{{\bf K}(s)(y-x)}{\bf C}_1^{-1}{\bf 1}\{y>x\}dy,\nonumber\\
{\bf E}(dy)^{(x)}(s)_{12}&=&{\bf\Psi}(s){\bf E}(dy)^{(x)}(s)_{22}
+e^{{\bf K}(s)(y-x)}{\bf\Psi}(s){\bf C}_2^{-1}{\bf 1}\{y> x\}dy,\nonumber
\\
{\bf E}(dy)^{(x)}(s)_{10}&=&
\left[
\begin{array}{cc}
{\bf E}(dy)^{(x)}(s)_{11}&{\bf E}(dy)^{(x)}(s)_{12}
\end{array}
\right]
\left[
\begin{array}{c}
{\bf T}_{10}\\
{\bf T}_{20}
\end{array}
\right]
\nonumber\\
&&
\times
(-({{\bf T}}_{00}-s{\bf I})^{-1})
dy.
\label{eq:Ed}
\end{eqnarray}
\end{lemma}
{\bf Proof:} The expressions for ${\bf E}(dy)^{(x)}(s)_{10}$ and ${\bf E}(dy)^{(x)}(s)_{20}$ follow by the argument in the proof of Lemma~\ref{lem:Edy_0}. Further, by partitioning the sample paths, since the process may visit level $y$ after returning to level $x$ first, or without hitting level $x$ at all, we have
\begin{eqnarray}
{\bf E}(dy)^{(x)}(s)_{11}&=&
{\bf\Psi}(s){\bf E}(dy)^{(x)}(s)_{21}
+{\bf E}(d(y-x))^{(0)}(s)_{11}{\bf 1}\{y>x\}dy,
\nonumber\\
{\bf E}(dy)^{(x)}(s)_{12}&=&{\bf\Psi}(s){\bf E}(dy)^{(x)}(s)_{22}
+{\bf E}(d(y-x))^{(0)}(s)_{12}{\bf 1}\{y>x\}dy
,\nonumber\\
\end{eqnarray}
and so the expressions for ${\bf E}(dy)^{(x)}(s)_{11}$ and ${\bf E}(dy)^{(x)}(s)_{12}$ follow by Lemma~\ref{lem:Edy_0}.

Next, we consider ${\bf E}(dy)^{(x)}(s)_{22}$. For $x>0$ define matrix ${\bf G}(x,t)=[G(x,t)_{kj}]$ such that for $k,j\in\mathcal{S}$,
\begin{equation}
G(x,t)_{kj}=P(\theta(0)\leq t,\varphi(\theta(0))=j\ |\
X(0)=x,\varphi(0)=k)
\end{equation}
is the probability that given the process starts from level $x$ in phase $k$, the process first hits level $0$ by time $t$, and does so in phase $j$. We partition ${\bf G}(x,t)$ according to $\mathcal{S}_1\cup\mathcal{S}_2$ as
\begin{equation}
{\bf G}(x,t)=
\left[
\begin{array}{cc}
{\bf 0}&{\bf G}(x,t)_{12}\\
{\bf 0}&{\bf G}(x,t)_{22}
\end{array}
\right].
\end{equation}
Also, define $\widetilde{\bf G}(x,s)=\int_{t=0}^{\infty}
e^{-s t}d{\bf G}(x,t)$, which we partition in an analogous manner.

The expression for ${\bf E}(dy)^{(x)}(s)_{22}$ then follows by partitioning the sample paths. The process can visit level $y$ in some phase in $\mathcal{S}_2$ directly after a visit to level $y$ in some phase in $\mathcal{S}_1$, or without visiting level $y$ in some phase in $\mathcal{S}_1$ at all, and so we take the sum of expressions corresponding to these two possibilities, which gives
\begin{eqnarray}
{\bf E}(dy)^{(x)}(s)_{22}&=&{\bf E}(dy)^{(x)}(s)_{21}{\bf C}_1{\bf\Psi}(s){\bf C}_2^{-1}
+\widetilde{\bf G}(x-y,s){\bf C}_2^{-1}{\bf 1}\{y<x\}dy,
\nonumber\\
\end{eqnarray}
and the result follows since by~\cite{BOT}, $\widetilde{\bf G}(x-y,s)=e^{{\bf D}(s)(x-y)}$.

Finally, we consider ${\bf E}(dy)^{(x)}(s)_{21}$. Denote
\begin{equation}
\underline{X}(t)=\inf_{u\in[0,t]}\{X(u)\}.
\end{equation}
Note that, given the process starts with $X(0)=x$, $\varphi(0)=i$, for the process to end with $X(t)\in dy$, $\varphi(t)=j$, with a taboo $\theta(0)>t$, one of the following two alternatives must occur.

The first alternative is that $y\geq x$. In this case,
\begin{itemize}
\item first, given $X(0)=x$, $\varphi(0)=i$, the process must reach some infimum $\underline{X}(t)=z\in(0,x]$ at some time $u\in[0,t]$, in some phase in $\mathcal{S}_2$, with the corresponding density recorded by matrix ${\bf G}_{22}(x-z,u)$; which is followed by an instantaneous transition to some phase $k$ in $\mathcal{S}_1$ according to the rate recorded by the block matrix ${\bf Q}_{21}$ of the fluid generator ${\bf Q}$, by the physical interpretation of ${\bf Q}$ in~\cite{BOT}. The corresponding density of this occurring is therefore $[{\bf G}_{22}(x-z,u){\bf Q}_{21}]_{ik}$.
\item Next, starting from level $z$ in phase $k$ at time $u$, the process must remain above level $z$ during the time interval $[u,t]$, ending in some level $y$ in phase $j$ at time $t$. The corresponding density of this occurring is $[\bphi(y-z,t-u)]_{kj}$.
\end{itemize}

Consequently, the LST of this alternative is
\begin{eqnarray*}
\lefteqn{
\int_{z=0}^x\int_{t=0}^{\infty}\int_{u=0}^t
e^{-s t}{\bf G}_{22}(x-z,u){\bf Q}_{21}(s)\bphi(y-z,t-u)_{11}dudtdz
}\\
&=&
\int_{z=0}^x\int_{u=0}^{\infty}\int_{t=u}^{\infty}
e^{-s u}
{\bf G}_{22}(x-z,u){\bf Q}_{21}(s)e^{-s (t-u)}\bphi(y-z,t-u)_{11}dudtdz\\
&=&
\int_{z=0}^x
\left(
\int_{u=0}^{\infty}e^{-s u}{\bf G}_{22}(x-z,u)du
\right)
{\bf Q}_{21}(s)
\left(
\int_{t=0}^{\infty}
e^{-s t}\bphi(y-z,t)_{11}dt
\right)
dz\\
&=&
\int_{z=0}^x
e^{{\bf D}(s)(x-z)}{\bf Q}_{21}(s)e^{{\bf K}(s)(y-z)}dz
.
\end{eqnarray*}

The second alternative is that $y<x$. The LST of this alternative, by an argument similar to above, is
\begin{eqnarray*}
\int_{z=0}^y
e^{{\bf D}(s)(x-z)}{\bf Q}_{21}(s)e^{{\bf K}(s)(y-z)}dz.
\end{eqnarray*}

Taking the sum of the expressions corresponding to the two alternatives and right-multiplying by ${\bf C}_1^{-1}$ results in the integral expression for ${\bf E}(dy)^{(x)}(s)_{21}$. \rule{9pt}{9pt}

\begin{remark}\label{howtogetEx21s} Consider
	\begin{equation}
	{\bf E}^{(x)}(s)_{21}
	=\int_{y=0}^{\infty}{\bf E}(dy)^{(x)}(s)_{21}={\bf X}{\bf C}_1^{-1},
	\end{equation}
where ${\bf X}=\int_{y=0}^{\infty}{\bf X}(y)dy$, and
\begin{eqnarray}\label{eq:Xy}
	{\bf X}(y)&=&
	\int_{z=0}^{\min\{x,y\}}
	e^{{\bf D}(s)(x-z)}{\bf Q}_{21}(s)e^{{\bf K}(s)(y-z)}dz.
	\end{eqnarray}
Then, by integration by parts in~\eqref{eq:Xy},
	${\bf X}(y)$ is the solution of
	\begin{eqnarray}\label{eq:Xyb}
	{\bf D}(s){\bf X}(y)+{\bf X}(y){\bf K}(s)
	&=&-\left[
	e^{{\bf D}(s)(x-z)}{\bf Q}_{21}(s)e^{{\bf K}(s)(y-z)}
	\right]_{z=0}^{\min\{x,y\}}.
	\end{eqnarray}
and by integrating~\eqref{eq:Xyb}, ${\bf X}$ is the solution of
	\begin{eqnarray}
	{\bf D}(s){\bf X}+{\bf X}{\bf K}(s)
&=&
e^{{\bf D}(s)(x)}{\bf Q}_{21}(s)(-{\bf K}(s)^{-1})
-(-{\bf D}(s)^{-1}){\bf Q}_{21}(s)
\nonumber\\
&&
+(-{\bf D}(s)^{-1})e^{{\bf D}(s)(x)}{\bf Q}_{21}(s)
+{\bf Q}_{21}(s){\bf K}(s)^{-1}
.
	\end{eqnarray}	
	
\end{remark}

\section{Approach}\label{sec:heaviside}
\newcommand{\imi}{{\rm i}}

The key idea is to write each of ${\bf E}(dy)^{(x)}(s)$ and ${\bf E}^{(x)}(s)$ in the form
\begin{equation}\label{eq:form}
\tilde f(s)=\tilde f(s^*) - C (s-s^*)^{1/2}+o((s-s^*)^{1/2}),
\end{equation}
and then apply the {\em Heaviside principle} in order to evaluate~(\ref{mu_fraction}). In this section, we summarise the relevant mathematical background required for this analysis.

Consider a function $f:\mathbb{R}\to\mathbb{R}$.
Let $\tilde f(s):=\int_0^\infty e^{-s x}f(x)\ d
x$ for $s\in \mathbb{R}$ be its Laplace transform. Consider singularities of $\tilde{f}(s)$.
We assume that one with the largest strictly negative real
part is real and we denote it by $s^*<0$. Notice that this yields the integrability of
$\int_0^\infty|f(x)|\ dx$. The inversion formula reads
\begin{equation}
f(x)=\frac{1}{2\pi \imi}\int_{a-\imi \infty}^{a+\imi\infty}\tilde{f}(s) e^{sx}\ ds
\end{equation}
for some (and then any) $a>s^*$.

We now focus on a class of theorems that infer the tail behaviour of a function from its Laplace transform, commonly referred to as {\it Tauberian theorems}. Importantly, the behaviour of the Laplace transform around the singularity $s^*$ plays a crucial role here. The following heuristic  principle given in \cite{AW1997} is often relied upon. Suppose that for $s^*$, some constants $K$ and $C$,  and a non-integer $q > 0$,
\begin{equation}
\tilde{f}(s)=K - C
(s-s^*)^q+o((s-s^*)^q),\qquad
\mbox{as $s\downarrow s^*$.}
\end{equation}
Then
\begin{equation}
f(x)=\frac{C}{\Gamma(-q)} x^{-q-1}e^{s^*
	x}(1+o(1)),\qquad \mbox{as $x\to\infty$,}
\end{equation}
where $\Gamma(\cdot)$ is the gamma function. Below we specify conditions under which this relation can be rigorously proven. Later in our paper we apply it for the specific case that $q=1/2$; recall that $\Gamma(-1/2)=-2\sqrt{\pi}$.

A formal justification of the above relation can be found in
Doetsch~\cite[Theorem~37.1]{Doetsch1974}. Following Miyazawa and
Rolski~\cite{miyazawaRolski}, we consider the following specific
form. For this we  first recall the concept of the $\mathfrak
W$-contour with an half-angle of opening $\pi/2<\psi\le \pi$, as
depicted in~\cite[Fig. 30, p. 240]{Doetsch1974}; also, ${\mathscr
	G}_{\alpha}(\psi)$  is the region between the contour
$\mathfrak W$ and the line $\Re(z)=0$. More precisely,
\begin{equation}
{\mathscr G}_{\alpha}(\psi) \equiv \{z \in \mathbb{C}; \Re(z) < 0, z \ne \alpha,
\ |\ \arg (z - \alpha)| < \psi \},
\end{equation}
where $\arg z$ is the principal part of the argument of the complex number~$z$.
In the following theorem, conditions  are identified such that the
above principle holds; we refer to this as the {\it Heaviside's
	operational principle}, or simply {\it Heaviside principle}.

\begin{theorem}[Heaviside principle]\label{t.tauberian}
	Suppose that for  $\tilde f:\mathbb{C}\to \mathbb{C}$ and
	$s^* < 0$ the following three conditions hold:
	\begin{itemize}
		\item [{\em (A1)}] $\tilde f(\cdot)$ is analytic in a
		region ${\mathscr G}_{s^*}(\psi)$ for some $\pi/2<\psi\le
		\pi$;
		\item [{\em (A2)}] $\tilde f(s) \to 0$ as $|s| \to \infty$ with $s \in {\mathscr G}_{s^*}(\psi)$;
		\item [{\em (A3)}]  for some constants $K$ and $C$, and a non-integer $q>0$,
		\begin{eqnarray}
		\label{eqn:bridge case 1} \tilde f(s)=K - C (s-s^*)^{q}+o((s-s^*)^{q}),
		\end{eqnarray}
		where ${\mathscr G}_{s^*}(\psi)\ni s \to s^*$.
	\end{itemize}
	Then
	\begin{eqnarray*}
		f(x)=\frac{C}{\Gamma(-q)} x^{-q-1}e^{s^* x}(1+o(1))
	\end{eqnarray*}
	as $x\to\infty$.
\end{theorem}

We now discuss when assumption (A1) is satisfied.  To check that
the Laplace transform $\tilde{f}(\cdot)$ is  analytic in the
region $\mathscr G_{s^*}(\psi)$, we can use the concept of
semiexponentiality  of $f$  (see \cite[p. 314]{henrici}).
\begin{definition}[Semiexponentiality]
	$f$  is said to be {\em semiexponential} if for some  $0< \phi\le\pi/2$
and all $-\phi\le {\vartheta}\le \phi$
there exists finite
	and strictly negative
	$\gamma({\vartheta})$, defined  as the infimum of all such $a$ such that
	\[\left|f(e^{\imi{\vartheta}}r)\right|<e^{ar}\]
	for all sufficiently large $r$.
\end{definition}

Relying on this concept, the following sufficient condition for (A1) applies.

\begin{proposition}\label{a1semexp}{\em \cite[Thm. 10.9f]{henrici}}
	Suppose that $f$ is  semiexponential with $\gamma({\vartheta})$
	fulfilling the following conditions: {\em (i)}
	$\gamma=\gamma(0)<0$, {\em (ii)} $\gamma({\vartheta})\geq\gamma(0)$ in a neighborhood of ${\vartheta}=0$, and
	{\em (iii)} it is smooth. Then {\em (A1)} is satisfied.
\end{proposition}

Note that by Lemma \ref{th:Edy},
all assumptions of Proposition  \ref{a1semexp} are satisfied
and we can apply the Heaviside principle given in Theorem \ref{t.tauberian} for ${\bf E}(dy)^{(x)}(s)$ and ${\bf E}^{(x)}(s)$.

\section{Application of the Heaviside principle}
\label{sec:applH}

By Section~\ref{sec:LSTs}, ${\bf E}(dy)^{(x)}(s)$ and ${\bf E}^{(x)}(s)$ are expressed in terms of ${\bf Q}(s)$ and ${\bf\Psi}(s)$, and so we
derive the expansion around $s^*$ for each of them first.

Consider ${\bf\Psi}(s)$ defined in~(\ref{def_Psi_s}). We have ${\bf\Psi}(s)=\int_{t=0}^{\infty}e^{-st}\bpsi(t)dt<\infty$ for all $s\geq 0$ by~\cite{BOT,BOT5}. Define the singularity
\begin{eqnarray}
s^*&=&\max\{s\leq 0:{\bf\Psi}(s)<\infty, {\bf\Psi}(z)=\infty \mbox{ for all }z< s\},
\end{eqnarray}
where the existence of $s^*$ follows from \cite[Thm. 3.3, p. 15]{Doetsch1974}.

Consider matrices ${\bf K}(s)$ and ${\bf D}(s)$ defined in~(\ref{eq:KsDs}), and recall that $sp({\bf K}(s))\cap sp(-{\bf D}(s))=\varnothing$ for all $s\geq 0$. Define
\begin{eqnarray}\label{def:deltastar}
\delta^*&=&\max\{s\in[s^*,0):sp({\bf K}(s))\cap sp(-{\bf D}(s))\not=\varnothing\},
\end{eqnarray}
whenever the maximum exists. The definition implies that ${\bf K}(\delta^*)$ and $(-{\bf D}(\delta^*))$ have a common eigenvalue.

\begin{lemma}\label{lem:sstardeltastar2}
	We have $s^*=\delta^*$.
\end{lemma}
{\bf Proof:} Consider equation~(\ref{Ric}) and for all $s$ for which ${\bf Q}(s)$ exists, define function of ${\bf X}=[x_{ij}]_{i\in\mathcal{S}_1,j\in\mathcal{S}_2}$,
	\begin{equation}\label{Ric_ab}
	g_s({\bf X})={\bf Q}_{12}(s)+{\bf Q}_{11}(s){\bf X}+{\bf X}{\bf Q}_{22}(s)+{\bf X}{\bf Q}_{21}(s){\bf X},
	\end{equation}
	where, for ${\bf X}\geq {\bf 0}$, ${\bf X}\not= {\bf 0}$, we have
	\begin{eqnarray}\label{neg_der}
	\frac{d}{ds}g_s({\bf X})&=&
	-{\bf C}^{-1}_1{\bf T}_{10}({\bf T}_{00}-s{\bf I})^{-2}{\bf T}_{02}
	{\bf X}
	-{\bf C}^{-1}_1\left({\bf I}+{\bf T}_{10}({\bf T}_{00}-s{\bf I})^{-2}{\bf T}_{01}\right)
	{\bf X}
	\nonumber\\
	&&
	-{\bf X}
	{\bf C}^{-1}_2\left({\bf I}+{\bf T}_{10}({\bf T}_{00}-s{\bf I})^{-2}{\bf T}_{02}\right)
	-
	{\bf X}
	{\bf C}^{-1}_2{\bf T}_{10}({\bf T}_{00}-s{\bf I})^{-2}{\bf T}_{01}
	{\bf X}
	\nonumber\\
	&<&{\bf 0},
	\end{eqnarray}
	since $({\bf T}_{00}-s{\bf I})^{-2}=\left(\int_{t=0}^{\infty}e^{-st}e^{{\bf T}_{00}t}dy\right)^2>{\bf 0}$, and so $g_s({\bf X})$ is a decreasing function of $s$.
	
	Also, define functions $g_s^{(u,v)}({\bf X})=[g_s({\bf X})]_{uv}$, for $u\in\mathcal{S}_1$, $v\in\mathcal{S}_2$,
	\begin{eqnarray}\label{Ric_ab_surface}
	g_s^{(u,v)}({\bf X})&=&[{\bf Q}_{12}(s)]_{u,v}
	+\sum_{k\in\mathcal{S}_1} [{\bf Q}_{11}(s)]_{uk}x_{kv}
	+\sum_{\ell\in\mathcal{S}_2} x_{u\ell}[{\bf Q}_{22}(s)]_{\ell v}
	\nonumber\\
	&&
	+\sum_{k\in\mathcal{S}_2,\ell\in\mathcal{S}_1} x_{uk}[{\bf Q}_{21}(s)]_{k\ell}x_{\ell v},
	\end{eqnarray}
	each corresponding to an $|\mathcal{S}_1|\times|\mathcal{S}_2|$-dimensional quadratic
	smooth surface. The matrix equation~(\ref{Ric}) is equivalent to the system of $|\mathcal{S}_1|\times|\mathcal{S}_2|$ quadratic polynomial equations, given by,
	\begin{eqnarray}\label{Ric_ab_surface2}
	g_s^{(u,v)}({\bf X})&=&0 \mbox{ for all }u\in\mathcal{S}_1,v\in\mathcal{S}_2,
	\end{eqnarray}
	each corresponding to the $(u,v)$-th level curve.
	
	Now, by Lemma~\ref{elm:psi}, for all $s\geq s^*$, ${\bf\Psi}(s)$ is a solution of $g_{s}({\bf X})={\bf 0}$ and so is an intersection point of all level curves~\eqref{Ric_ab_surface2}.
		
	Some other solutions to $g_{s}({\bf X})={\bf 0}$ may exist. For all real $s$, we denote by ${\bf X}(s)$ the family of solutions that correspond to the intersection point ${\bf\Psi}(s)$. That is, when $s\geq s^*$, ${\bf X}(s)={\bf\Psi}(s)$, and if ${\bf X}(s)$ exists for $s<s^*$ in some neighbourhood of $s^*$, then ${\bf X}(s)$ must be a continuous function of $s$ in such neighbourhood, due to the monotonicity and continuity of $g_s({\bf X})$.
		
	So suppose that there exist solutions ${\bf X}(s)$ to $g_{s}({\bf X})={\bf 0}$ for $s<s^*$ in some neighbourhood of $s^*$, and that $\lim_{s\uparrow s^*}{\bf X}(s)={\bf\Psi}(s^*)$. Then, since ${\bf\Psi}(s^*)>{\bf 0}$, there exists ${\bf W}>{\bf 0}$ with $g_{s}({\bf W})={\bf 0}$ for some $s<s^*$ with $sp({\bf Q}_{11}(s))\cap sp(-{\bf Q}_{22}(s))=\varnothing$ (due to the fact that spectra $sp({\bf Q}_{11}(s))$ and $sp(-{\bf Q}_{22}(s))$ are discrete).

	Therefore, by~\cite[Theorem 2.3]{Guo2002} and~\cite[Algorithm 1]{BOT5}, we have ${\bf\Psi}(s)<\infty$ for such $s<s^*$, and this contradicts the definition of $s^*$. Consequently, ${\bf X}(s)$ does not exist for $s<s^*$, and so the level curves~(\ref{Ric_ab_surface2}) must touch (have a common tangent line) at $s=s^*$, but not at $s>s^*$.

	Denote
	\begin{equation}
		\nabla g_{s^*}^{(u,v)}(x_{ij},[x_{ij}^*])
			=
			\frac{\partial}{\partial x_{ij}}g^{(u,v)}([x_{ij}])\Big|_{[x_{ij}^*]},
	\end{equation}	
	and note that
		\begin{eqnarray}
			\lefteqn{
			\nabla g_{s^*}^{(u,v)}(x_{ij},[x_{ij}^*])
			=
			\frac{\partial}{\partial x_{ij}}g_{s^*}^{(u,v)}([x_{ij}])\Big|_{[x_{ij}^*]}
		}
		\nonumber
			\\
			&=&
			[{\bf Q}_{11}(s^*)]_{ui}{\bf 1}\{j=v\}
			+[{\bf Q}_{22}(s^*)]_{jv}{\bf 1}\{i=u\}
			+\sum_{k\in\mathcal{S}_2} [{\bf Q}_{21}(s^*)]_{ki}\times x_{uk}^*   {\bf 1}\{j=v\}
			\nonumber\\
			&&+
			\sum_{\ell\in\mathcal{S}_1}[{\bf Q}_{21}(s^*)]_{j\ell}\times x_{\ell v}^*{\bf 1}\{i=u\}
			\nonumber\\
			&=&
			[{\bf K}(s^*)]_{ui}{\bf 1}\{j=v\}
			+[{\bf D}(s^*)]_{jv}{\bf 1}\{i=u\}
			.\nonumber\\\label{eq:grad}
			\end{eqnarray}

	The tangent plane to the $(u,v)$-th level curve~(\ref{Ric_ab_surface2}) at ${\bf\Psi}(s^*)=[x_{ij}^*]$, is the solution to the equation,
		\begin{eqnarray}
		0&=&
		\sum_{i,j} \nabla g^{(u,v)}_{s^*}(x_{ij},[x_{ij}^*]) (x_{ij}-x_{ij}^*)
		\nonumber\\
		&=&\sum_{i,j}
				\Big(
				[{\bf K}(s^*)]_{ui}{\bf 1}\{j=v\}
						+[{\bf D}(s^*)]_{jv}{\bf 1}\{i=u\}
						\Big) (x_{ij}-x_{ij}^*)
						\nonumber\\
					&=&
		\sum_i[{\bf K}(s^*)]_{ui}[{\bf X}-{\bf\Psi}(s^*)]_{iv}
		+\sum_j [{\bf X}-{\bf\Psi}(s^*)]_{uj}[{\bf D}(s^*)]_{jv}
		\nonumber\\
		&=&
		\left[  {\bf K}(s^*)({\bf X}-{\bf\Psi}(s^*))
		+({\bf X}-{\bf\Psi}(s^*)){\bf D}(s^*)\right]_{uv}.
		\label{eq:tang_keep}
		\end{eqnarray}

	From linear algebra, a matrix equation of the form ${\bf 0}={\bf A}{\bf X}+{\bf X}{\bf B}$ has a nonzero solution if and only if ${\bf A}$ and $(-{\bf B})$ have a common eigenvalue (e.g. see~\cite{Bha-Ros}). Therefore, the equation
	\begin{eqnarray}\label{xx}
	{\bf 0}&=&
	{\bf K}(s^*)({\bf X}-{\bf\Psi}(s^*))
	+({\bf X}-{\bf\Psi}(s^*)){\bf D}(s^*)
	\end{eqnarray}
	has a solution ${\bf Z}=[z_{ij}]\not= {\bf\Psi}(s^*)$ if and only if ${\bf K}(s^*)$ and $(-{\bf D}(s^*))$ have a common eigenvalue, in which case the tangent planes~\eqref{eq:tang_keep} to all level curves~(\ref{Ric_ab_surface2}) at ${\bf\Psi}(s^*)$, intersect with one another at a tangent line that goes through ${\bf Z}$ and ${\bf\Psi}(s^*)$.
	
	That is, the level curves~(\ref{Ric_ab_surface2}) touch if and only if $sp({\bf K}(s^*))\cap sp(-{\bf D}(s^*))\not=\varnothing$. Hence, $s^*=\delta^*$. \rule{9pt}{9pt}

We now extend the result for $s>0$ in~\cite[Theorem 1]{BOT} to all $s\geq s^*$.
\begin{cor} For all $s\geq s^*$, ${\bf\Psi}(s)$ is the minimum nonnegative solution of the Riccati equation~\eqref{Ric}. 
\end{cor}
{\bf Proof:}\  Suppose $s\geq s^*$. Then, ${\bf Q}_{11}(s)\leq {\bf K}(s)={\bf Q}_{11}(s)+{\bf\Psi}(s){\bf Q}_{21}(s)$ and ${\bf Q}_{22}(s)\leq {\bf D}(s)={\bf Q}_{22}(s)+{\bf Q}_{21}(s){\bf\Psi}(s)$, and so $sp({\bf Q}_{11}(s))\cap sp(-{\bf Q}_{22}(s))=\varnothing$.

Therefore, by~\cite[Theorem 2.3]{Guo2002} and~\cite[Algorithm 1]{BOT5}, ${\bf\Psi}(s)$ is the minimum nonnegative solution of~\eqref{Ric}. \rule{9pt}{9pt}

In order to illustrate the theory, we consider the following simple example, which we will analyse as we develop the results throughout the paper.
\begin{example}\label{ex1} \rm
	Let $\mathcal{S}=\{1,2\}$, $\mathcal{S}_1=\{1\}$, $\mathcal{S}_2=\{2\}$, $c_1=1$, $c_2=-1$, and
	\begin{eqnarray}
	{\bf T}
	&=&\left[
	\begin{array}{cc}
	{\bf T}_{11}&{\bf T}_{12}\\
	{\bf T}_{21}&{\bf T}_{22}
	\end{array}
	\right]
	=\left[
	\begin{array}{cc}
	-a&a\\
	b&-b
	\end{array}
	\right],
	\\
	{\bf Q}(s)
	&=&\left[
	\begin{array}{cc}
	{\bf Q}_{11}(s)&{\bf Q}_{12}(s)\\
	{\bf Q}_{21}(s)&{\bf Q}_{22}(s)
	\end{array}
	\right]
	=\left[
	\begin{array}{cc}
	-a-s&a\\
	b&-b-s
	\end{array}
	\right],
	\end{eqnarray}
	with $a>b>0$ so that the process is stable.

	Then ${\bf\Psi}(s)$ is the minimum nonnegative solution of~(\ref{Ric}), here equivalent to
	\begin{equation}
	b x^2 - (a+b+2s)x+a=0,
	\end{equation}
	which  has solutions provided $\Delta(s)=(a+b+2s)^2-4ab \geq 0$, that is, for all
	\begin{eqnarray}
	s\in
	\left(-\infty,
	\frac{-(a+b)-2\sqrt{ab}}{2}
	\right]
	\cup\left[
	\frac{-(a+b)+2\sqrt{ab}}{2},
	+\infty\right)
	.
	\end{eqnarray}

	Since
	\begin{eqnarray}
	(a+b+2s)-\sqrt{\Delta(s)}\geq 0
	\iff
	s\leq \frac{2ab}{a+b},
	\end{eqnarray}
	it follows that ${\bf\Psi}(s)$ exists for all $s\geq \frac{-(a+b)+2\sqrt{ab}}{2}$, and
	\begin{eqnarray}
	{\bf\Psi}(s)&=&\frac{(a+b+2s)-\sqrt{\Delta(s)}}{2b},\\
	{\bf K}(s)&=&-a-s+\frac{(a+b+2s)-\sqrt{\Delta(s)}}{2},\\
	{\bf D}(s)&=&-b-s+\frac{(a+b+2s)-\sqrt{\Delta(s)}}{2}.
	\end{eqnarray}
	
	Therefore,
	\begin{equation}
	s^*=\frac{-(a+b)+2\sqrt{ab}}{2}<0
	\label{eq:s_star_ex1}
	\end{equation}
	and
	\begin{eqnarray}
	{\bf\Psi}(s^*)&=&\sqrt{\frac{a}{b}},\\
	{\bf K}(s^*)&=&-a-s^*+\sqrt{\frac{a}{b}}\ b
	=
	\frac{b-a}{2}<0
	,
	\\
	{\bf D}(s^*)&=&-b-s^*+b\ \sqrt{\frac{a}{b}}
	=\frac{a-b}{2}>0,
	\end{eqnarray}
	and note that $s^*=\delta^*$.
	
\end{example}

\begin{lemma}\label{QS_form}
For all $s>s^*$,
\begin{eqnarray}
{\bf Q}_{22}(s)&=&{\bf Q}_{22}(s^*)-{\bf A}_{22}(s^*)(s-s^*)+o(s-s^*),\\
{\bf Q}_{11}(s)&=&{\bf Q}_{11}(s^*)-{\bf A}_{11}(s^*)(s-s^*)+o(s-s^*),\\
{\bf Q}_{12}(s)&=&{\bf Q}_{12}(s^*)-{\bf A}_{12}(s^*)(s-s^*)+o(s-s^*),\\
{\bf Q}_{21}(s)&=&{\bf Q}_{21}(s^*)-{\bf A}_{21}(s^*)(s-s^*)+o(s-s^*),
\end{eqnarray}
where, for all $s> s^*$,
\begin{eqnarray}\label{eq:Aasder}
{\bf A}_{22}(s)&=&
-\frac{d}{ds}{\bf Q}_{22}(s)
=
{\bf C}^{-1}_2\left({\bf I}+{\bf T}_{20}({\bf T}_{00}-s{\bf I})^{-2}{\bf T}_{02} \right),\\
{\bf A}_{11}(s)&=&
-\frac{d}{ds}{\bf Q}_{11}(s)
=
{\bf C}^{-1}_1\left({\bf I}+{\bf T}_{10}({\bf T}_{00}-s{\bf I})^{-2}{\bf T}_{01} \right),\\
{\bf A}_{12}(s)&=&
-\frac{d}{ds}{\bf Q}_{12}(s)
=
{\bf C}^{-1}_1{\bf T}_{10}({\bf T}_{00}-s{\bf I})^{-2}{\bf T}_{02},\\
{\bf A}_{21}(s)&=&
-\frac{d}{ds}{\bf Q}_{21}(s)
=
{\bf C}^{-1}_2{\bf T}_{20}({\bf T}_{00}-s{\bf I})^{-2}{\bf T}_{01},
\end{eqnarray}
and ${\bf A}_{22}(s^*)=\lim_{s\downarrow s^*}{\bf A}_{22}(s)<\infty$, ${\bf A}_{11}(s^*)=\lim_{s\downarrow s^*}{\bf A}_{11}(s)<\infty$, ${\bf A}_{12}(s^*)=\lim_{s\downarrow s^*}{\bf A}_{12}(s)<\infty$, and ${\bf A}_{21}(s^*)=\lim_{s\downarrow s^*}{\bf A}_{21}(s)<\infty$.
\end{lemma}
{\bf Proof:}\ For all $s>s^*$,
	\begin{eqnarray}	
	-\frac{d}{ds}\left(-({\bf T}_{00}-s{\bf I})^{-1} \right)
	&=&
({\bf T}_{00}-s{\bf I})^{-2},
	\end{eqnarray}
and so by~(\ref{eq:Qs}),
\begin{eqnarray}
-\frac{d}{ds}{\bf Q}_{22}(s)
&=&
{\bf C}^{-1}_2\left({\bf I}+{\bf T}_{20}({\bf T}_{00}-s{\bf I})^{-2}{\bf T}_{02} \right),
\end{eqnarray}
which, with the notation ${\bf A}_{22}(s)=-\frac{d}{ds}{\bf Q}_{22}(s)$, implies,
\begin{eqnarray}\label{eqQ22}
{\bf Q}_{22}(s+h)={\bf Q}_{22}(s)-{\bf A}_{22}(s)(h)+o(h).
\end{eqnarray}

Next, since ${\bf\Psi}(s^*)<\infty$ by Lemma~\ref{lem:sstardeltastar2}, we have $\left(-({\bf T}_{00}-s^*{\bf I})^{-1} \right)<\infty$ and $({\bf T}_{00}-s^*{\bf I})^{-2}<\infty$, which implies that ${\bf A}_{22}(s^*)<\infty$. Taking the limits as $s\downarrow s^*$ in~(\ref{eqQ22}), and substituting $h=(s-s^*)$ gives
\begin{eqnarray}\label{eqA22}
{\bf Q}_{22}(s)&=&{\bf Q}_{22}(s^*)-{\bf A}_{22}(s^*)(s-s^*)+o(s-s^*).
\end{eqnarray}
The proof of the remaining expressions is analogous. \rule{9pt}{9pt}

For $s>s^*$, let
\begin{equation}\label{defPhi}
{\bf \Phi } (s)=\frac{d}{ds}{\bf \Psi}(s)=\lim_{h\to 0}\frac{{\bf \Psi}(s+h)-{\bf \Psi}(s)}{h},
\end{equation}
and, for $s\geq s^*$, let
\begin{eqnarray}\label{Us*}
{\bf U}(s)&=&
{\bf A}_{12}(s)+{\bf A}_{11}(s){\bf\Psi}(s)
+{\bf\Psi}(s){\bf A}_{22}(s)
+{\bf\Psi}(s){\bf A}_{21}(s){\bf\Psi}(s),
\nonumber\\
\end{eqnarray}
noting that ${\bf U}(s^*)$ exists by Lemma~\ref{QS_form}.
\begin{lemma}\label{cor_nosol}
	For $s>s^*$, ${\bf \Phi }(s)$ is the unique solution of the equation
	\begin{eqnarray}\label{eq_nosol}
	{\bf K}(s){\bf X}+{\bf X}{\bf D}(s)
	&=&
	{\bf U}(s).
	\end{eqnarray}
Furthermore, ${\bf \Phi }(s^*)=\lim_{s\downarrow s^*}{\bf \Phi }(s)=-\infty$.
\end{lemma}	
{\bf Proof:} By Lemma~\ref{lem:sstardeltastar2}, for all $s>s^*$, ${\bf K}(s)$ and $(-{\bf D}(s))$ have no common eigenvalues, and so by~\cite[Theorem 13.18]{laub}, the equation~(\ref{eq_nosol}) has a unique solution. We now show that ${\bf \Phi } (s)$ is the solution of~(\ref{eq_nosol}). Also see~\cite[Corollary 3]{BOT}. Indeed, by taking derivatives w.r.t. $s$ in the equation~(\ref{Ric}) for ${\bf \Psi}(s)$, we have
\begin{eqnarray}\label{eg:Phi}
{\bf 0}&=&
\frac{d}{ds}
\Big(
{\bf Q}_{12}(s)+{\bf Q}_{11}(s){\bf \Psi}(s)+{\bf \Psi}(s){\bf Q}_{22}(s)+{\bf \Psi}(s){\bf Q}_{21}(s){\bf \Psi}(s)
\Big)
\nonumber\\
&=&
-{\bf A}_{12}(s)
-{\bf A}_{11}(s){\bf \Psi}(s)+{\bf Q}_{11}(s){\bf \Phi}(s)
-{\bf \Psi}(s){\bf A}_{22}(s)+{\bf \Phi}(s){\bf Q}_{22}(s)
\nonumber\\
&&
+{\bf \Phi}(s){\bf Q}_{21}(s){\bf \Psi}(s)
-{\bf \Psi}(s){\bf A}_{21}(s){\bf \Psi}(s)
+{\bf \Psi}(s){\bf Q}_{21}(s){\bf \Phi}(s)
\nonumber\\
&=&
-{\bf U}(s)+{\bf K}(s){\bf \Phi } (s)+{\bf \Phi } (s){\bf D}(s).
\end{eqnarray}
Also, ${\bf \Phi}(s)<{\bf 0}$, since
\begin{equation}
{\bf \Phi}(s)=\frac{d}{ds}{\bf\Psi}(s)=\frac{d}{ds}\int_0^{\infty}e^{-st}\bpsi(t)dt=-\int_0^{\infty}te^{-st}\bpsi(t)dt<{\bf 0}.
\end{equation}

When $s=s^*$ however, by Lemma~\ref{lem:sstardeltastar2}, ${\bf K}(s)$ and $(-{\bf D}(s))$ have a common eigenvalue, and so by~\cite[Theorem 13.18]{laub}, the equation~(\ref{eq_nosol}) does not have a unique solution.

Finally, we show that $\lim_{s\downarrow s^*}{\bf \Phi }(s)=-\infty$. By standard methodology~\cite[Section 13.3]{laub}, for $s>s^*$, the unique solution to the equation~(\ref{eq_nosol}) can be written in the form
	\begin{eqnarray}\label{eq:uniquePhis^*}
	vec({\bf \Phi }(s))
	&=&
	({\bf Z}(s))^{-1}vec({\bf U}(s))
	=
	\frac{adj({\bf Z}(s))}{det({\bf Z}(s))}vec({\bf U}(s)),
	\end{eqnarray}
	where $vec({\bf \Phi }(s))$ and $vec({\bf U}(s))$ are column vectors obtained by stacking the columns (from the left to the right) of the original matrices one under another,
	\begin{equation}
	{\bf Z}(s)=({\bf I}\otimes {\bf K}(s))+({\bf D}(s)^T \otimes {\bf I}),
	\end{equation}
	and the eigenvalues of ${\bf Z}(s)$ are $(\lambda_i-\mu_j)$, where $\lambda_i$ are eigenvalues of ${\bf K}(s)$ and $\mu_j$ are eigenvalues of $(-{\bf D}(s))$. Since $det({\bf Z}(s))$ is the product of the eigenvalues of ${\bf Z}(s)$, and as $s\downarrow s^*$ one of the eigenvalues will approach zero due to $s^*=\delta^*$ by Lemma~\ref{lem:sstardeltastar2}, we have $\lim_{s\downarrow s^*}det({\bf Z}(s))=0$  and so ${\bf \Phi }(s^*)=\lim_{s\downarrow s^*}{\bf \Phi }(s)=-\infty$, where the negative sign is due to ${\bf \Phi }(s)<{\bf 0}$ for all $s>s^*$. \rule{9pt}{9pt}

We now state the key result of this paper.
\begin{theorem}\label{lem:psi_form}
	For all $s>s^*$,
	\begin{eqnarray}\label{eq:eqB}
		{\bf\Psi}(s)&=&{\bf\Psi}(s^*)-{\bf B}(s^*)\sqrt{s-s^*}+o(\sqrt{s-s^*}),
	\end{eqnarray}
	where ${\bf 0}<{\bf B}(s^*)<\infty$ solves
	\begin{eqnarray}
		{\bf B}(s^*){\bf Q}_{21}(s^*){\bf B}(s^*) &=&{\bf U}(s^*)-{\bf Y}(s^*),\label{eq1}\\
		{\bf K}(s^*){\bf B}(s^*)+{\bf B}(s^*){\bf D}(s^*)&=&{\bf 0},
		\label{eq2}
	\end{eqnarray}
	and
	\begin{eqnarray}
	{\bf Y}(s^*)&=&\lim_{s\downarrow s^*}
	\left(
	{\bf K}(s^*){\bf \Phi}(s)+{\bf \Phi}(s){\bf D}(s^*)
	\right)
	\label{whatisY}.
	\end{eqnarray}
\end{theorem}
{\bf Proof:}\ Note that for any function $h(\cdot)$ with $h(s-s^*)=o(s-s^*)$ or $h(s-s^*)=c\cdot (s-s^*)$ for some constant $c$, we have
\begin{equation}
	-\lim_{s\downarrow s^*}
	\left(
	\frac{{\bf \Psi}(s)-{\bf \Psi}(s^*)}{s-s^*}h(s-s^*)
	\right)
	={\bf 0}.
\end{equation}
Consider $h(s-s^*)=(s-s^*)/|| {\bf \Psi}(s)-{\bf \Psi}(s^*) ||$. We have,
\begin{equation}
	\lim_{s\downarrow s^*}\frac{s-s^*}{h(s-s^*)}
	=\lim_{s\downarrow s^*}|| {\bf \Psi}(s)-{\bf \Psi}(s^*) ||
	=0,
\end{equation}
which implies $(s-s^*)=o(h(s-s^*))$, and
\begin{equation}
	\lim_{s\downarrow s^*}
	\Bigg|\Bigg|
	\frac{{\bf \Psi}(s)-{\bf \Psi}(s^*)}{s-s^*}h(s-s^*)
	\Bigg|\Bigg|
	=1\not= 0.
\end{equation}
Therefore, there exists a continuous, positive-valued function $h(\cdot)$ such that $(s-s^*)=o(h(s-s^*))$ and
\begin{equation}
	-\lim_{s\downarrow s^*}
	\left(
	\frac{{\bf \Psi}(s)-{\bf \Psi}(s^*)}{s-s^*}h(s-s^*)
	\right)
	={\bf B}(s^*)
\end{equation}
for some constant matrix ${\bf 0}<{\bf B}(s^*)<\infty$. For such $h(\cdot)$, define function $g(\cdot)$ such that
\begin{equation}
	g(s-s^*)=\frac{s-s^*}{h(s-s^*)},
\end{equation}
with clearly $\lim_{s\downarrow s^*}g(s-s^*)=0$ since $(s-s^*)=o(h(s-s^*))$.

Consequently, we have
\begin{equation}
	-\lim_{s\downarrow s^*}
	\left(
	\frac{{\bf \Psi}(s)-{\bf \Psi}(s^*)}{g(s-s^*)}
	\right)
	={\bf B}(s^*),
\end{equation}
which implies that
\begin{eqnarray}
	{\bf\Psi}(s)&=&{\bf\Psi}(s^*)-{\bf B}(s^*)g(s-s^*)
	+o(g(s-s^*)).
\end{eqnarray}

We now solve for ${\bf B}(s^*)$ and $g(s-s^*)$. By~(\ref{Ric}) and Lemma~\ref{QS_form}, since
\begin{eqnarray}
	{\bf 0}&=&
	{\bf Q}_{12}(s^*)
	+ {\bf Q}_{11}(s^*){\bf\Psi}(s^*)
	+{\bf\Psi}(s^*){\bf Q}_{22}(s^*)
	+{\bf\Psi}(s^*){\bf Q}_{21}(s^*){\bf\Psi}(s^*),
\end{eqnarray}
we have,
\begin{eqnarray}
{\bf 0}
&=&
	{\bf Q}_{12}(s)
	+{\bf Q}_{11}(s){\bf\Psi}(s)+{\bf\Psi}(s){\bf Q}_{22}(s)
	+{\bf\Psi}(s){\bf Q}_{21}(s){\bf\Psi}(s)
\nonumber
	\\
	&=&\left({\bf Q}_{12}(s^*)-{\bf A}_{12}(s^*)(s-s^*)\right)
	\nonumber\\
	&&
	+
	\left( {\bf Q}_{11}(s^*)-{\bf A}_{11}(s^*)(s-s^*) \right)
	\left( {\bf\Psi}(s^*)-{\bf B}(s^*)g(s-s^*) \right)
	\nonumber\\
	&&
	+ \left( {\bf\Psi}(s^*)-{\bf B}(s^*)g(s-s^*) \right)
	\left( {\bf Q}_{22}(s^*)-{\bf A}_{22}(s^*)(s-s^*) \right)
	\nonumber\\
	&&
	+
	\left( {\bf\Psi}(s^*) - {\bf B}(s^*)g(s-s^*) \right)
	\left( {\bf Q}_{21}(s^*)-{\bf A}_{21}(s^*)(s-s^*) \right)
	\left( {\bf\Psi}(s^*) - {\bf B}(s^*)g(s-s^*) \right)
	\nonumber\\
	&&
	+o(s-s^*)+o(g(s-s^*)),\nonumber\\
\end{eqnarray}
and so
\begin{eqnarray}\label{keyeq}
	{\bf 0}
	&=&
	-(s-s^*){\bf U}(s^*)+g(s-s^*){\bf W}(s^*)
	+g^2(s-s^*){\bf V}(s^*)
	+o(s-s^*)
	+o(g(s-s^*)),
	\nonumber\\
\end{eqnarray}
where ${\bf U}(s^*)$ is defined in~(\ref{Us*}), and
\begin{eqnarray}\label{eq:WV}
	{\bf W}(s^*)&=&
	\left(
	{\bf Q}_{11}(s^*)+{\bf\Psi}(s^*){\bf Q}_{21}(s^*)
	\right)
	{\bf B}(s^*)
	+{\bf B}(s^*)
	\left(
	{\bf Q}_{22}(s^*)+{\bf Q}_{21}(s^*){\bf\Psi}(s^*)
	\right)
	\nonumber\\
	&=&
	{\bf K}(s^*){\bf B}(s^*)+{\bf B}(s^*){\bf D}(s^*)
	,
	\nonumber\\
	{\bf V}(s^*)
	&=&
	{\bf B}(s^*){\bf Q}_{21}(s^*){\bf B}(s^*)
	.\nonumber\\
\end{eqnarray}

We now use equation~(\ref{keyeq}) in order to solve for ${\bf B}(s^*)$ and $g(s-s^*)$. We note that ${\bf V}(s^*)\not= {\bf 0}$ and ${\bf U}(s^*)\not= {\bf 0}$. Indeed, ${\bf V}(s^*)\not= {\bf 0}$ since ${\bf V}(s^*)> {\bf 0}$ due to ${\bf B}(s^*)>{\bf 0}$, ${\bf Q}_{21}(s^*)\geq {\bf 0}$, ${\bf Q}_{21}(s^*)\not= {\bf 0}$. Further, ${\bf U}(s^*)\not= {\bf 0}$ since ${\bf U}(s^*)> {\bf 0}$. Indeed, in the case $\mathcal{S}_0=\varnothing$, since ${\bf C}_1,{\bf C}_2>{\bf 0}$ and ${\bf\Psi}(s^*)>{\bf 0}$, we have ${\bf U}(s^*)={\bf C}_1^{-1}{\bf \Psi}(s^*)+{\bf \Psi}(s^*){\bf C}_2^{-1}>{\bf 0}$. In the case $\mathcal{S}_0\not=\varnothing$, we have $-({\bf T}_{00}-s^*{\bf I})^{-1}=\int_{t=0}^{\infty}e^{-s^*t}e^{{\bf T}_{00}t}dt>{\bf 0}$, and $({\bf T}_{00}-s^*{\bf I})^{-2}=(-({\bf T}_{00}-s^*{\bf I})^{-1})^2>0$. Therefore ${\bf A}_{11}(s^*),{\bf A}_{22}(s^*)>{\bf 0}$, ${\bf A}_{12}(s^*),{\bf A}_{21}(s^*)\geq 0$ and ${\bf\Psi}(s^*)>{\bf 0}$, and so ${\bf U}(s^*)={\bf A}_{12}(s^*)+{\bf A}_{11}(s^*){\bf\Psi}(s^*)
+{\bf\Psi}(s^*){\bf A}_{22}(s^*)
+{\bf\Psi}(s^*){\bf A}_{21}(s^*){\bf\Psi}(s^*)>{\bf 0}$. 

Consequently, below we consider two cases, ${\bf W}(s^*)\not= {\bf 0}$ and ${\bf W}(s^*)= {\bf 0}$, respectively, labelled Case I and Case II below.

\bigskip
{\bf Case I.}\ Suppose ${\bf W}(s^*)\not= {\bf 0}$. Then,
\begin{eqnarray}\label{eq:caseI}
	{\bf 0}
	&=&
	-(s-s^*){\bf U}(s^*)+g(s-s^*){\bf W}(s^*)
	+g^2(s-s^*){\bf V}(s^*)
	+o(s-s^*)
	+o(g(s-s^*)).
	\nonumber\\
\end{eqnarray}

Consider $(s-s^*)$ and $g(s-s^*)$. Either one of them dominates another, or one is a multiple of the other.

(i)\ If $g(s-s^*)=o(s-s^*)$, then dividing equation~(\ref{eq:caseI}) by $(s-s^*)$ and taking limits as $s\downarrow s^*$ gives ${\bf 0}={\bf U}(s^*)$, a contradiction.

(ii)\ If $(s-s^*)=o(g(s-s^*))$, then dividing equation~(\ref{eq:caseI}) by $g(s-s^*)$ and taking limits as $s\downarrow s^*$ gives ${\bf 0}={\bf W}(s^*)$, a contradiction.

(iii)\ If $g(s-s^*)=c\cdot (s-s^*)$ for some constant $c>0$, then without loss of generality we may assume $c=1$, since ${\bf B}(s^*)g(s-s^*)=({\bf B}(s^*)c )(s-s^*)$ suggests the substitution $\tilde {\bf B}(s^*)\equiv {\bf B}(s^*)c$. Then we have,
\begin{eqnarray}\label{eq:Bs}
	{\bf\Psi}(s)&=&{\bf\Psi}(s^*)-\tilde{\bf B}(s^*)(s-s^*)
	+o(s-s^*),
\end{eqnarray}
with $\tilde {\bf B}(s^*)<\infty$. However, dividing equation~(\ref{eq:Bs}) by $(s-s^*)$ and taking limits as $s\downarrow s^*$ gives, by Lemma~\ref{cor_nosol},
\begin{eqnarray}
	\tilde{\bf B}(s^*)=-\lim_{s\downarrow s^*}\frac{{\bf \Psi}(s)-{\bf \Psi}(s^*)}{s-s^*}=\infty,
\end{eqnarray}
a contradiction.

That is, the assumption ${\bf W}(s^*)\not= {\bf 0}$ leads to a contradiction.

\bigskip
{\bf Case II.}\ By above, we must have ${\bf W}(s^*)= {\bf 0}$, or equivalently,
\begin{equation}
{\bf K}(s^*){\bf B}(s^*)+{\bf B}(s^*){\bf D}(s^*)={\bf 0},
\end{equation}
and so,
\begin{eqnarray}\label{eqII}
	{\bf 0}&=&-(s-s^*){\bf U}(s^*)+g^2(s-s^*){\bf V}(s^*)+o(s-s^*)
	+o(g(s-s^*)).
\end{eqnarray}
We note that $g^2(s-s^*)=o(g(s-s^*))$, and consider the following.

{\bf(i)}\ First, we show that $(s-s^*)=o(g(s-s^*))$. Indeed, if $g(s-s^*)=o(s-s^*)$ or $g(s-s^*)=c\cdot (s-s^*)$ for some $c\not= 0$, then dividing equation~(\ref{eqII}) by $(s-s^*)$ and taking limits as $s\downarrow s^*$ gives ${\bf U}(s^*)={\bf 0}$, a contradiction. Therefore we must have $(s-s^*)=o(g(s-s^*))$. That is, $g(s-s^*)$ dominates both $g^2(s-s^*)$ and $(s-s^*)$.

Then,
\begin{eqnarray}
	\lim_{s\to s^*}\frac{o(s-s^*)}{g(s-s^*)}=\lim_{s\to s^*}\frac{o(s-s^*)}{(s-s^*)}\frac{(s-s^*)}{g(s-s^*)}
	=0,
\end{eqnarray}
which gives $o(s-s^*)=o(g(s-s^*))$, and so we write~(\ref{eqII}) in the form
\begin{eqnarray}\label{eqIImore}
	{\bf 0}&=&-(s-s^*){\bf U}(s^*)+g^2(s-s^*){\bf V}(s^*)
	+o(g(s-s^*)).
\end{eqnarray}

Since $(s-s^*)=o(g(s-s^*))$, we consider two cases, $o(g(s-s^*))={\bf 0}$ and $o(g(s-s^*))\not ={\bf 0}$, respectively, labelled (A) and (B) below.

{\bf(A)}\ Suppose $o(g(s-s^*))={\bf 0}$. Then~(\ref{eqIImore}) reduces to
\begin{eqnarray}\label{eqIImore_zero}
{\bf 0}&=&-(s-s^*){\bf U}(s^*)+g^2(s-s^*){\bf V}(s^*).
\end{eqnarray}
If $(s-s^*)=o(g^2(s-s^*))$, we divide~(\ref{eqIImore_zero}) by $g^2(s-s^*)$ and take limits as $s\downarrow s^*$ to get ${\bf V}(s^*)={\bf 0}$, a contradiction. If $g^2(s-s^*)=o(s-s^*)$, we divide~(\ref{eqIImore_zero}) by $(s-s^*)$ and take limits as $s\downarrow s^*$ to get ${\bf U}(s^*)={\bf 0}$, a contradiction. So we must have $(s-s^*)=c\cdot g^2(s-s^*)$ for some constant $c>0$, and without loss of generality we may assume $c=1$. Then, dividing equation~(\ref{eqIImore_zero}) by $(s-s^*)$ and taking limits as $s\downarrow s^*$ gives ${\bf U}(s^*)={\bf V}(s^*)$, or equivalently,
\begin{equation}
{\bf B}(s^*){\bf Q}_{21}(s^*){\bf B}(s^*)={\bf U}(s^*).
\end{equation}
That is, Case (A) gives $g(s-s^*)=\sqrt{s-s^*}$.

{\bf(B)}\ Suppose $o(g(s-s^*))\not={\bf 0}$. Then we write the term $o(g(s-s^*))$ in the form
\begin{eqnarray}
	o(g(s-s^*))&=& L(s-s^*){\bf Y}(s^*) +o(L(s-s^*))
\end{eqnarray}
for some function $L(\cdot)\not= 0$ such that $L(s-s^*)=o(g(s-s^*))$ and some constant ${\bf Y}(s^*)\not= {\bf 0}$.

Then we have,
\begin{eqnarray}\label{eqIImorev2}
	{\bf 0}&=&-(s-s^*){\bf U}(s^*)+g^2(s-s^*){\bf V}(s^*)
	+L(s-s^*){\bf Y}(s^*) +o(L(s-s^*)).
	\nonumber\\
\end{eqnarray}
Consider the terms $(s-s^*)$, $g^2(s-s^*)$ and $L(s-s^*)$, and the following cases under assumption~(B), labelled (B)(ii)-(B)(iv), respectively. We will show that Case (B)(ii) gives a contradiction and Cases (B)(iii)-(iv) give $g(s-s^*)=\sqrt{s-s^*}$. 

{\bf(B)(ii)}\  Suppose one of $(s-s^*)$, $g^2(s-s^*)$ and $L(s-s^*)$, dominates the two others.

If $(s-s^*)$ dominates the two others, that is $g^2(s-s^*)=o(s-s^*)$ and $L(s-s^*)=o(s-s^*)$, then dividing equation~(\ref{eqIImorev2}) by $(s-s^*)$ and taking limits as $s\downarrow s^*$ gives ${\bf U}(s^*)={\bf 0}$, a contradiction.

If $g^2(s-s^*)$ dominates the two others, that is $(s-s^*)=o(g^2(s-s^*))$ and $L(s-s^*)=o(g^2(s-s^*))$, then dividing equation~(\ref{eqIImorev2}) by $g^2(s-s^*)$ and taking limits as $s\downarrow s^*$ gives ${\bf V}(s^*)={\bf 0}$, a contradiction.

If $L(s-s^*)$ dominates the two others, that is $(s-s^*)=o(L(s-s^*))$ and $g^2(s-s^*)=o(L(s-s^*))$, then dividing equation~(\ref{eqIImorev2}) by $L(s-s^*)$ and taking limits as $s\downarrow s^*$ gives ${\bf Y}(s^*)={\bf 0}$, a contradiction.

That is, Case (B)(ii) gives a contradiction. Therefore at least two of $(s-s^*)$, $g^2(s-s^*)$ and $L(s-s^*)$ must be a multiple of each other.

{\bf(B)(iii)}\ Suppose each of $(s-s^*)$, $g^2(s-s^*)$ and $L(s-s^*)$ is a multiple of any other. Then, $(s-s^*)=c\cdot g^2(s-s^*)=d\cdot L(s-s^*)$, and without loss of generality we may assume $c=1$, $d=1$, by argument analogous to before. Therefore, dividing equation~(\ref{eqIImorev2}) by $(s-s^*)$ and taking limits as $s\downarrow s^*$ gives ${\bf 0}=-{\bf U}(s^*)+{\bf V}(s^*)+{\bf Y}(s^*)$, or equivalently,
\begin{eqnarray}
{\bf B}(s^*){\bf Q}_{21}(s^*){\bf B}(s^*) &=&{\bf U}(s^*)-{\bf Y}(s^*).
\end{eqnarray}
That is, Case (B)(iii) gives $g(s-s^*)=\sqrt{s-s^*}$. 

{\bf(B)(iv)}\ Suppose exactly two of $(s-s^*)$, $g^2(s-s^*)$ and $L(s-s^*)$ are a multiple of one another. Then such two terms must dominate the third term, or we have a contradiction by part (i) of Case II above.

If $(s-s^*)=c\cdot g^2(s-s^*)$ for some $c>0$, then without loss of generality we may assume $c=1$. Also, we must have $L(s-s^*)=o(s-s^*)$. Therefore, $g(s-s^*)=\sqrt{s-s^*}$, and dividing equation~(\ref{eqIImorev2}) by $(s-s^*)$ and taking limits as $s\downarrow s^*$ gives ${\bf U}(s^*)={\bf V}(s^*)$, or equivalently,
\begin{equation}
{\bf B}(s^*){\bf Q}_{21}(s^*){\bf B}(s^*)={\bf U}(s^*).
\end{equation}

If $L(s-s^*)=c\cdot (s-s^*)$ for some $c\not= 0$, then dividing equation~(\ref{eqIImorev2}) by $(s-s^*)$ and taking limits as $s\downarrow s^*$ gives ${\bf V}(s^*)={\bf 0}$, a contradiction.

If $L(s-s^*)=c\cdot g^2(s-s^*)$ for some $c> 0$, then without loss of generality we may assume $c=1$. Also, we must have $L(s-s^*)=o(s-s^*)$. Therefore, equation~(\ref{eqIImorev2}) becomes,
\begin{eqnarray}\label{eqIImorev3}
	{\bf 0}&=&-(s-s^*){\bf U}(s^*)+g^2(s-s^*)({\bf V}(s^*)+{\bf Y}(s^*))
	+o(g^2(s-s^*)).
\end{eqnarray}
In this case, if $g^2(s-s^*)=o(s-s^*)$ then dividing equation~(\ref{eqIImorev3}) by $(s-s^*)$ and taking limits as $s\downarrow s^*$  gives ${\bf U}(s^*)={\bf 0}$, a contradiction. If $(s-s^*)=o(g^2(s-s^*))$ then dividing equation~(\ref{eqIImorev3}) by $g^2(s-s^*)$ and taking limits as $s\downarrow s^*$  gives ${\bf U}(s^*)+{\bf Y}(s^*)={\bf 0}$, a contradiction. Therefore, we must have $g^2(s-s^*)=c\cdot (s-s^*)$ for some $c> 0$. Without loss of generality we may assume $c=1$. Therefore, $g(s-s^*)=\sqrt{s-s^*}$, and dividing equation~(\ref{eqIImorev3}) by $g^2(s-s^*)$ and taking limits as $s\downarrow s^*$  gives
\begin{eqnarray}
	{\bf B}(s^*){\bf Q}_{21}(s^*){\bf B}(s^*)&=&
	{\bf U}(s^*)- {\bf Y}(s^*).
\end{eqnarray}
That is, Case (B)(iv) gives $g(s-s^*)=\sqrt{s-s^*}$. 

By above cases, we must have
\begin{equation}
	g(s-s^*)=\sqrt{s-s^*},
\end{equation}
and
\begin{eqnarray}
	{\bf B}(s^*){\bf Q}_{21}(s^*){\bf B}(s^*)&=&
	{\bf U}(s^*)- {\bf Y}(s^*),
	\\
	{\bf K}(s^*){\bf B}(s^*)+{\bf B}(s^*){\bf D}(s^*)&=&{\bf 0},
\end{eqnarray}
and $-\infty<{\bf Y}(s^*)\leq {\bf U}(s^*)$. Here, ${\bf Y}(s^*)= {\bf 0}$ whenever the term $o(g(s-s^*))$ in~(\ref{eqIImore}) satisfies $o(g(s-s^*))=o(s-s^*)$, and ${\bf Y}(s^*)\not= {\bf 0}$ when $o(g(s-s^*))=(s-s^*){\bf Y}(s^*)+o(s-s^*)$.

Finally, we show~\eqref{whatisY}. By L'Hospital's rule,
\begin{eqnarray}
{\bf B}(s^*)&=&
-\lim_{s\downarrow s^*}
\frac{{\bf \Psi}(s)-{\bf \Psi}(s^*)}
{\sqrt{s-s^*}}
\ = \
-\lim_{s\downarrow s^*}
\left(
{\bf \Phi}(s)2\sqrt{s-s^*}\
\right)
,
\end{eqnarray}
and so, by taking limits as $s\downarrow s^*$ in~(\ref{eg:Phi}),
\begin{eqnarray}
\lefteqn{
{\bf B}(s^*){\bf Q}_{21} (s^*){\bf B}(s^*)
+{\bf Y}(s^*)
=
	{\bf U}(s^*)
}
\nonumber\\
&=&
	\lim_{s\downarrow s^*}
	\Big[
	{\bf K}(s){\bf \Phi}(s)
	+{\bf \Phi}(s){\bf D}(s)
	\Big]
\nonumber\\
&=&
\lim_{s\downarrow s^*}
\Big[
({\bf Q}_{11}(s)+{\bf \Psi}(s){\bf Q}_{21}(s)){\bf \Phi}(s)
+{\bf \Phi}(s)({\bf Q}_{22}(s)+{\bf Q}_{21}(s){\bf \Psi}(s))
\Big]
\nonumber\\
&=&
\lim_{s\downarrow s^*}
\Bigg[
\frac{1}{2}
\left(
\frac{{\bf \Psi}(s)-{\bf \Psi}(s^*)}{\sqrt{s-s^*}}
\right)
{\bf Q}_{21}(s)
\left(
{\bf \Phi}(s)2\sqrt{s-s^*}
\right)
\nonumber\\
&&
+
\frac{1}{2}
\left(
{\bf \Phi}(s)2\sqrt{s-s^*}
\right)
{\bf Q}_{21}(s)
\left(
\frac{{\bf \Psi}(s)-{\bf \Psi}(s^*)}{\sqrt{s-s^*}}
\right)
\nonumber\\
&&
+
{\bf Q}_{11}(s){\bf \Phi}(s)+{\bf \Phi}(s){\bf Q}_{22}(s)
+{\bf \Psi}(s^*){\bf Q}_{21}(s){\bf \Phi}(s)
+{\bf \Phi}(s){\bf Q}_{21}(s){\bf \Psi}(s^*)
\Bigg]
\nonumber\\
&=&
\frac{1}{2}{\bf B}(s^*){\bf Q}_{21} (s^*){\bf B}(s^*)
+\frac{1}{2}{\bf B}(s^*){\bf Q}_{21} (s^*){\bf B}(s^*)
\nonumber\\
&&
+
\lim_{s\downarrow s^*}
\Bigg[
\left(\frac{{\bf Q}_{11}(s)-{\bf Q}_{11}(s^*)}{\sqrt{s-s^*}}\right)
({\bf \Phi}(s)\sqrt{s-s^*})
+(\sqrt{s-s^*}{\bf \Phi}(s))\left(\frac{{\bf Q}_{22}(s)-{\bf Q}_{22}(s^*)}{\sqrt{s-s^*}} \right)
\nonumber\\
&&
+
{\bf Q}_{11}(s^*){\bf \Phi}(s)+{\bf \Phi}(s){\bf Q}_{22}(s^*)
\nonumber\\
&&+{\bf \Psi}(s^*)
\left(\frac{{\bf Q}_{21}(s)-{\bf Q}_{21}(s^*)}{\sqrt{s-s^*}}\right)
({\bf \Phi}(s)\sqrt{s-s^*})
+({\bf \Phi}(s)\sqrt{s-s^*})
\left(\frac{{\bf Q}_{21}(s)-{\bf Q}_{21}(s^*)}{\sqrt{s-s^*}}\right)
{\bf \Psi}(s^*)
\nonumber\\
&&
+{\bf \Psi}(s^*){\bf Q}_{21}(s^*){\bf \Phi}(s)+{\bf \Phi}(s){\bf Q}_{21}(s^*){\bf \Psi}(s^*)
\Bigg]
\nonumber\\
&=&
{\bf B}(s^*){\bf Q}_{21} (s^*){\bf B}(s^*)
+ {\bf 0}+
\lim_{s\downarrow s^*}
\Big[
{\bf K}(s^*){\bf \Phi}(s)+{\bf \Phi}(s){\bf D}(s^*)
\Big]
,\nonumber\\
\end{eqnarray}
which completes the proof. \rule{9pt}{9pt}

The next result follows immediately by Theorem~\ref{t.tauberian}.
\begin{corollary}\label{cr:psit}
We have
\begin{equation}
\bpsi(t)
={\bf B}(s^*)\frac{1}{2\sqrt{\pi}}t^{-3/2}e^{s^* t}(1+o(1)).
\end{equation}
\end{corollary}

\setcounter{example}{0}	
\begin{example}(continued)\label{ex1B} \rm
Since $\lim_{s\downarrow s^*} \Delta(s)=\Delta(s^*)=0$, we have
\begin{eqnarray}
\lim_{s\downarrow s^*}\frac{d}{ds}{\bf\Psi}(s)
&=& \lim_{s\downarrow s^*}\frac{d}{ds} \frac{(a+b+2s)-\sqrt{\Delta(s)}}{2b}
\nonumber\\
&=& \lim_{s\downarrow s^*}\frac{d}{ds}
\left(
\frac{1}{b}
-\frac{1}{4b \sqrt{\Delta(s)}}(8s+4(a+b))
\right)
\nonumber\\
&=&-\infty,
\end{eqnarray}
as expected. Furthermore,
\begin{eqnarray}
\lim_{s\downarrow (s^*)}\frac{{\bf\Psi}(s)-{\bf\Psi}(s^*)}{\sqrt{s-s^*}}
&=&
\lim_{s\downarrow (s^*)}
\left(
\frac{(a+b+2s)-\sqrt{\Delta(s)}}{2b\sqrt{s-s^*}}
-
\frac{(a+b+2s^*)-\sqrt{\Delta(s^*)}}{2b\sqrt{s-s^*}}
\right)
\nonumber\\
&=&	
\frac{\sqrt{2\sqrt{ab}}}{-b}
,
\end{eqnarray}
which implies
\begin{equation}
{\bf\Psi}(s)={\bf\Psi}(s^*)-{\bf B}(s^*)\sqrt{s-s^*}+o(\sqrt{s-s^*}),
\end{equation}
where
\begin{equation}
{\bf B}(s^*)=\frac{\sqrt{2\sqrt{ab}}}{b}.
\end{equation}

Therefore, by Theorem~\ref{t.tauberian},
\begin{equation}
\bpsi(t)
=
\frac{\sqrt{2\sqrt{ab}}}{2b\sqrt{\pi}}
t^{-3/2}
\exp\left(\left(
\frac{-(a+b)+2\sqrt{ab}}{2}
\right)
t
\right)
(1+o(1))
.
\end{equation}

Also, ${\bf A}_{12}(s^*)=0$, ${\bf A}_{21}(s^*)=0$, ${\bf A}_{11}(s^*)=1$, ${\bf A}_{22}(s^*)=1$, and so
\begin{eqnarray}
{\bf B}(s^*){\bf Q}_{21}(s^*){\bf B}(s^*)&=&
b\left(\frac{\sqrt{2\sqrt{ab}}}{b}\right)^2
= \ 2\sqrt{\frac{a}{b}}
\ = \ {\bf U}(s^*)
,
\end{eqnarray}
and
\begin{eqnarray}
{\bf K}(s^*){\bf B}(s^*)+{\bf B}(s^*){\bf D}(s^*)&=&
\left( \frac{b-a}{2}+\frac{a-b}{2}  \right)\left(\frac{\sqrt{2\sqrt{ab}}}{b}\right)
={\bf 0},
\\
	\lim_{s\downarrow s^*}
	\left(
	{\bf K}(s^*){\bf \Phi}(s)+{\bf \Phi}(s){\bf D}(s^*)
	\right)
&=&
\lim_{s\downarrow s^*}
\left(
\left( \frac{b-a}{2}+\frac{a-b}{2}  \right)
2 {\bf \Phi}(s)
\right)
={\bf 0}.
	\end{eqnarray}

\end{example}	

\medskip

Define matrices, for $n\geq 1$, 
\begin{eqnarray}
{\bf H}_{1,n}(s^*)&=&
\sum_{i=0}^{n-1}
\left({\bf K}(s^*) \right)^i
\times
{\bf B}(s^*){\bf Q}_{21}(s^*)
\times
\left({\bf K}(s^*) \right)^{n-1-i}
,
\end{eqnarray}
and
\begin{equation}
{\bf H}(s^*,y)=\sum_{n=1}^{\infty}
\frac{y^n}{n!}
{\bf H}_{1,n}(s^*), 
\ \
{\bf H}(s^*)=\int_{y=0}^{\infty}{\bf H}(s^*,y)dy,
\end{equation}
and a column vector
\begin{eqnarray}
\widetilde{\bf H}(s^*)&=&{\bf H}(s^*){\bf C}_1^{-1}{\bf 1}
+
\left(
-({\bf K}(s^*))^{-1}{\bf B}(s^*)
+
{\bf H}(s^*){\bf\Psi}(s^*)
\right)
{\bf C}_2^{-1}{\bf 1}
\nonumber\\
&&
+
\left[
\begin{array}{cc}
{\bf H}(s^*)&
-({\bf K}(s^*))^{-1}{\bf B}(s^*)
+
{\bf H}(s^*){\bf\Psi}(s^*)
\end{array}
\right]
\left[
\begin{array}{c}
{\bf C}_1^{-1}{\bf T}_{10}\\
{\bf C}_2^{-1}{\bf T}_{20}
\end{array}
\right]
(-({{\bf T}}_{00}-s^*{\bf I})^{-1}){\bf 1}.
\nonumber\\
\label{colvectildeHstar}
\end{eqnarray}

Below, we derive the expressions for $\bmu(dy)^{(0)}$. 
\begin{theorem}\label{thm:main}
The matrix $\bmu(dy)^{(0)}$ is unique and
\begin{eqnarray}
\bmu(dy)^{(0)}_{11}&=& 
diag({\widetilde{\bf H}(s^*)})^{-1}
{ {\bf H}(s^*,y){\bf C}_1^{-1}}  
dy,
\nonumber\\
\bmu(dy)^{(0)}_{12}&=& 
diag({\widetilde{\bf H}(s^*)})^{-1}
\left(
e^{{\bf K}(s^*)y}{\bf B}(s^*)
+
{\bf H}(s^*,y){\bf\Psi}(s^*)
\right){\bf C}_2^{-1}
dy,
\nonumber\\
\bmu(dy)^{(0)}_{10}&=&	
\left[
\begin{array}{cc}
\bmu(dy)^{(0)}_{11} &\bmu(dy)^{(0)}_{12}
\end{array}
\right]
\left[
\begin{array}{c}
{\bf T}_{10}\\
{\bf T}_{20}
\end{array}
\right]
(-({{\bf T}}_{00}-s^*{\bf I})^{-1}).
\nonumber\\
\end{eqnarray}	
\end{theorem}

\begin{remark}
From Theorem \ref{thm:main} it follows that
the crucial step in identifying Yaglom limit
given above is identification of $s^*$. Unfortunately, this must be done
for each stochastic fluid queue separately.
\end{remark}

{\bf Proof:}\ By Lemma~\ref{lem:Edy_0}, Lemma~\ref{QS_form} and Theorem~\ref{lem:psi_form}, we have
\begin{eqnarray}
\lefteqn{
e^{{\bf K}(s)y}
=
\lim_{K\to\infty}
\sum_{n=0}^{K}
\frac{y^n}{n!}
\left(
{\bf Q}_{11}(s)+{\bf\Psi}(s){\bf Q}_{21}(s)
\right)^n
}
\nonumber\\
&=&
\lim_{K\to\infty}\sum_{n=0}^{K}
\frac{y^n}{n!}
\left(
{\bf Q}_{11}(s^*)
+\left( {\bf\Psi}(s^*)-{\bf B}(s^*)\sqrt{s-s^*} \right)
{\bf Q}_{21}(s^*)
+o(\sqrt{s-s^*})
\right)^n
\nonumber\\
&=&
\lim_{K\to\infty}\sum_{n=0}^{K}
\frac{y^n}{n!}
\left(
{\bf Q}_{11}(s^*)+{\bf\Psi}(s^*){\bf Q}_{21}(s^*)
\right)^n
-\lim_{K\to\infty}\sqrt{s-s^*}\sum_{n=1}^{K }
\frac{y^n}{n!}
{\bf H}_{1,n}
+o(\sqrt{s-s^*})
\nonumber\\
&=&
e^{{\bf K}(s^*)y}
-\sqrt{s-s^*}{\bf H}(s^*,y)
+o(\sqrt{s-s^*}),
\nonumber\\
\label{edoK}
\end{eqnarray}
which gives
\begin{eqnarray}
{\bf E}(dy)^{(0)}(s)_{11}&=&
e^{{\bf K}(s)y}{\bf C}_1^{-1}dy
\nonumber\\
&=&{\bf E}(dy)^{(0)}(s^*)_{11}
-\sqrt{s-s^*}{\bf H}(s^*,y){\bf C}_1^{-1}dy
+o(\sqrt{s-s^*})
\label{Edyterms}
\end{eqnarray}
and
\begin{eqnarray}
\lefteqn{
{\bf E}(dy)^{(0)}(s)_{12}=
e^{{\bf K}(s)y}{\bf\Psi}(s){\bf C}_2^{-1}dy
}
\nonumber\\
&=&
\left(
e^{{\bf K}(s^*)y}
-\sqrt{s-s^*}{\bf H}(s^*,y)
\right)
({\bf\Psi}(s^*)-{\bf B}(s^*)\sqrt{s-s^*}){\bf C}_2^{-1}dy
+o(\sqrt{s-s^*})
\nonumber\\
&=&
{\bf E}(dy)^{(0)}(s^*)_{12}
-\sqrt{s-s^*}
\left(
e^{{\bf K}(s^*)y}{\bf B}(s^*)
+
{\bf H}(s^*,y){\bf\Psi}(s^*)
\right){\bf C}_2^{-1}dy
+o(\sqrt{s-s^*}).
\nonumber\\
\label{Edyterms2}
\end{eqnarray}
and, by noting that $({\bf T}_{00}-s{\bf I})^{-1}
-({\bf T}_{00}-s^*{\bf I})^{-1}
=
(s-s^*)({\bf T}_{00}-s^*{\bf I})^{-2}+o(s-s^*)$, which gives
$({{\bf T}}_{00}-s{\bf I})^{-1}=({{\bf T}}_{00}-s^*{\bf I})^{-1}+o(\sqrt{s-s^*})$, we have
\begin{eqnarray}
\lefteqn{
{\bf E}(dy)^{(0)}(s)_{10}
=
\left[
\begin{array}{cc}
	e^{{\bf K}(s)y}	{\bf C}_1^{-1}&e^{{\bf K}(s)y}{\bf\Psi}(s){\bf C}_2^{-1}
\end{array}
\right]
\left[
\begin{array}{c}
{\bf T}_{10}\\
	{\bf T}_{20}
\end{array}
\right]
(-({{\bf T}}_{00}-s{\bf I})^{-1})
}
\nonumber\\
&=&
{\bf E}(dy)^{(0)}(s^*)_{10}
-
\sqrt{s-s^*}
\left[
\begin{array}{cc}
{\bf H}(s^*,y)&
\left(
e^{{\bf K}(s^*)y}{\bf B}(s^*)
+
{\bf H}(s^*,y){\bf\Psi}(s^*)
\right)
\end{array}
\right]
\nonumber\\
&&\times
\left[
\begin{array}{c}
{\bf C}_1^{-1}{\bf T}_{10}\\
{\bf C}_2^{-1}{\bf T}_{20}
\end{array}
\right]
(-({{\bf T}}_{00}-s^*{\bf I})^{-1})
dy
+o(\sqrt{s-s^*}).
\nonumber\\
\label{Edyterms3}
\end{eqnarray}

Furthermore,
\begin{eqnarray}
{\bf E}^{(0)}(s)_1&=&
\int_{y=0}^{\infty}{\bf E}(dy)^{(0)}(s)_{11}{\bf 1}
+\int_{y=0}^{\infty}{\bf E}(dy)^{(0)}(s)_{12}{\bf 1}
+\int_{y=0}^{\infty}{\bf E}(dy)^{(0)}(s)_{10}{\bf 1}
\nonumber\\
&=&
{\bf E}^{(0)}(s^*)_1
-\sqrt{s-s^*}
{\bf H}(s^*){\bf C}_1^{-1}
{\bf 1}
-\sqrt{s-s^*}
\left(
-({\bf K}(s^*))^{-1}{\bf B}(s^*)
+
{\bf H}(s^*){\bf\Psi}(s^*)
\right){\bf 1}
\nonumber\\
&&
-\sqrt{s-s^*}
\left[
\begin{array}{cc}
{\bf H}(s^*)&
\left(
-({\bf K}(s^*))^{-1}{\bf B}(s^*)
+
{\bf H}(s^*){\bf\Psi}(s^*)
\right)
\end{array}
\right]
\nonumber\\
&&\times
\left[
\begin{array}{c}
{\bf C}_1^{-1}{\bf T}_{10}\\
{\bf C}_2^{-1}{\bf T}_{20}
\end{array}
\right]
(-({{\bf T}}_{00}-s^*{\bf I})^{-1}){\bf 1}
+o(\sqrt{s-s^*})
\nonumber\\
&=&
{\bf E}^{(0)}(s^*)_1
-\sqrt{s-s^*}
\widetilde{\bf H}(s^*)
+o(\sqrt{s-s^*}).
\nonumber\\
\label{E0s1}
\end{eqnarray}

The result follows by Theorem~\ref{t.tauberian} and \eqref{mu_fraction}, since the relevant terms cancel out. Indeed, for $i,j\in\mathcal{S}_1$, by~\eqref{Edyterms}-\eqref{E0s1}, Theorem~\ref{t.tauberian} and \eqref{mu_fraction},
	\begin{eqnarray}
	\mu(dy)^{(0)}_{ij}
	&=&
	\lim_{t\to\infty}
	\frac{
		P(X(t)\in dy,\varphi(t)=j,\theta(0)>t\ |\
		X(0)=0,\varphi(0)=i)
	}
	{
		P(\theta(0)>t\ |\ X(0)=0,\varphi(0)=i)
	}
	\nonumber\\
	&=&
	\frac{
		\lim_{t\to\infty}
		([{\bf H}(s^*,y)
		{\bf C}_1^{-1}]_{ij}\Gamma(1/2)^{-1}t^{-1/2-1}e^{s^*t}(1+o(1))
		}
		{
		\lim_{t\to\infty}
		([\widetilde{\bf H}(s^*)]_{i}\Gamma(1/2)^{-1}t^{-1/2-1}e^{s^*t}(1+o(1))
		}dy
		\nonumber\\
	&=&
		\frac{
[{\bf H}(s^*,y){\bf C}_1^{-1}]_{ij}
		}
		{
	[\widetilde{\bf H}(s^*)]_{i}
		}dy,
	\end{eqnarray}
	which gives the result for $\bmu(dy)^{(0)}_{11}$. Expressions for $\bmu(dy)^{(0)}_{12}$ and $\bmu(dy)^{(0)}_{10}$ follow in a similar manner. \rule{9pt}{9pt}

\setcounter{example}{0}
\begin{example}(continued)\label{ex1Bb} \rm Finally,
\begin{eqnarray}
{\bf H}(s^*,y)&=&
\sum_{n=1}^{\infty}
\frac{y^n}{n!}
{\bf H}_{1,n}(s^*)
\nonumber\\
&=&
\sum_{n=1}^{\infty}
\frac{y^n}{n!}
\sum_{i=0}^{n-1}
\left({\bf K}(s^*) \right)^i
\times
{\bf B}(s^*){\bf Q}_{21}(s^*)
\times
\left({\bf K}(s^*) \right)^{n-1-i}
\nonumber\\
&=&
\sum_{n=1}^{\infty}
\frac{y^n}{n!}
\sum_{i=0}^{n-1}
\left( -(a-b)/2\right)^{n-1}
\sqrt{2\sqrt{ab}}
\nonumber\\
&=&
ye^{\left( -(a-b)/2 \right)y}
\sqrt{2\sqrt{ab}},
\end{eqnarray}
and
\begin{eqnarray}
{\bf H}(s^*)&=&
\int_{y=0}^{\infty}
ye^{\left( -(a-b)/2 \right)y}
\sqrt{2\sqrt{ab}}\
dy
\nonumber\\
&=&
\frac{\sqrt{2\sqrt{ab}}}{(a-b)^2/4},
\nonumber\\
\widetilde{\bf H}(s^*)&=&{\bf H}(s^*)
+
\left(
-({\bf K}(s^*))^{-1}{\bf B}(s^*)
+
{\bf H}(s^*){\bf\Psi}(s^*)
\right)
\nonumber\\
&=&
\frac{\sqrt{2\sqrt{ab}}}{(a-b)^2/4}\
\left(1+\sqrt{\frac{a}{b}} \right)
+\frac{2}{a-b}
\left(
\frac{\sqrt{2\sqrt{ab}}}{b}
\right)
,
\end{eqnarray}
and so
\begin{eqnarray}
\bmu(dy)^{(0)}_{11}&=& 
	diag({\widetilde{\bf H}(s^*)})^{-1}
	{ {\bf H}(s^*,y)}  
	dy
	\nonumber\\
	&=&
\left( 
\frac{1}{(a-b)^2/4}\
	\left(1+\sqrt{\frac{a}{b}} \right)
	+\frac{2}{a-b}
	\left(
	\frac{1}{b}
	\right)
\right)^{-1}
ye^{(-(a-b)/2)y}
dy
,\nonumber\\
\bmu(dy)^{(0)}_{12}&=& 
diag({\widetilde{\bf H}(s^*)})^{-1}
\left(
e^{{\bf K}(s^*)y}{\bf B}(s^*)
+
{\bf H}(s^*,y){\bf\Psi}(s^*)
\right)
dy
\nonumber\\
&=&
	\left( 
\frac{1}{(a-b)^2/4}\
		\left(1+\sqrt{\frac{a}{b}} \right)
		+\frac{2}{a-b}
		\left(
		\frac{1}{b}
		\right)
	\right)^{-1}
	\left(
	\frac{1}{b}+y\sqrt{\frac{a}{b}}
	\right)	
	e^{(-(a-b)/2)y}
	dy
.
\nonumber\\
\end{eqnarray}
We plot the values of $\bmu(dy)^{(0)}_{11}$ and $\bmu(dy)^{(0)}_{12}$ in Figure~\ref{figex1}.
\begin{figure}[htbp]
\centering
\includegraphics[scale=0.5]{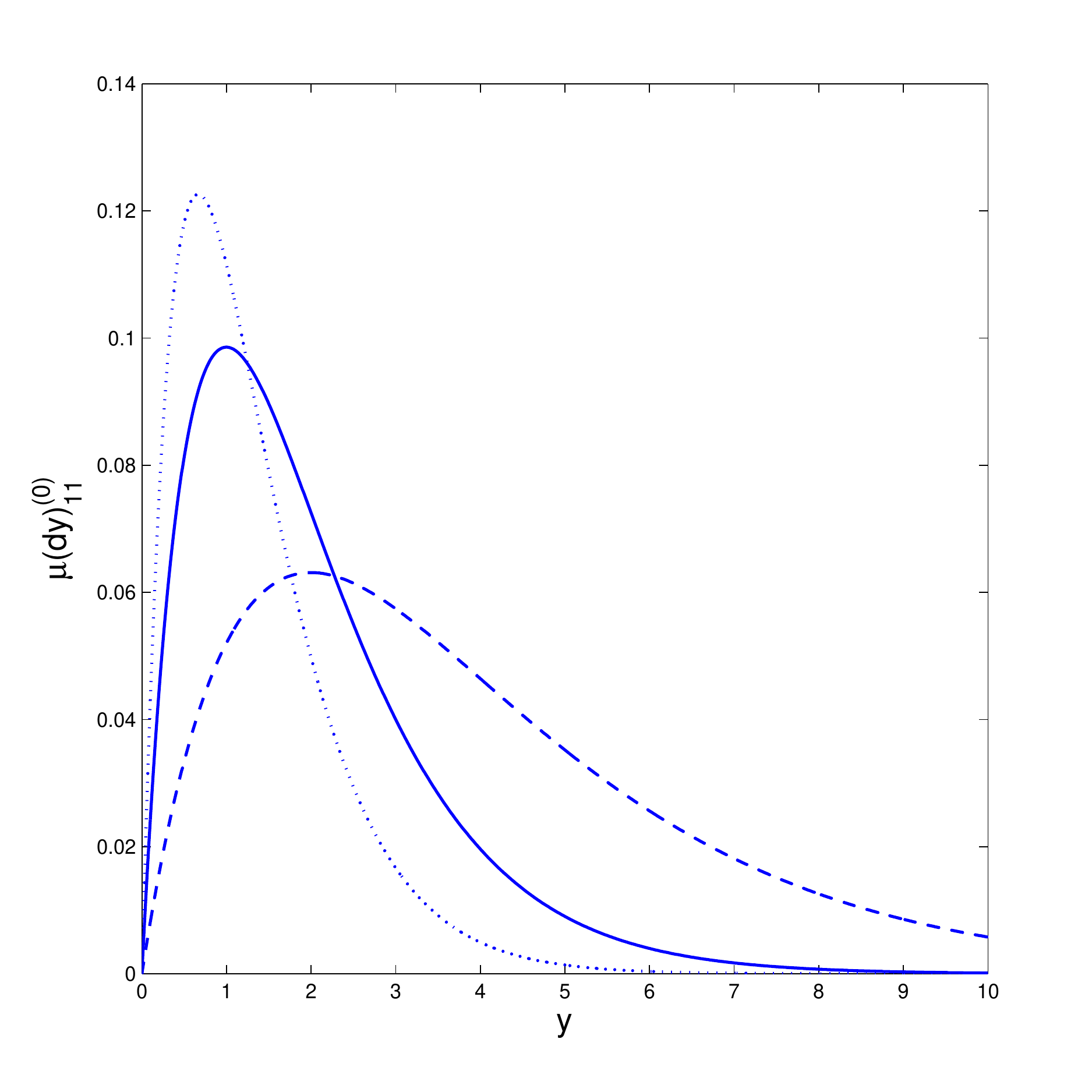}
\includegraphics[scale=0.5]{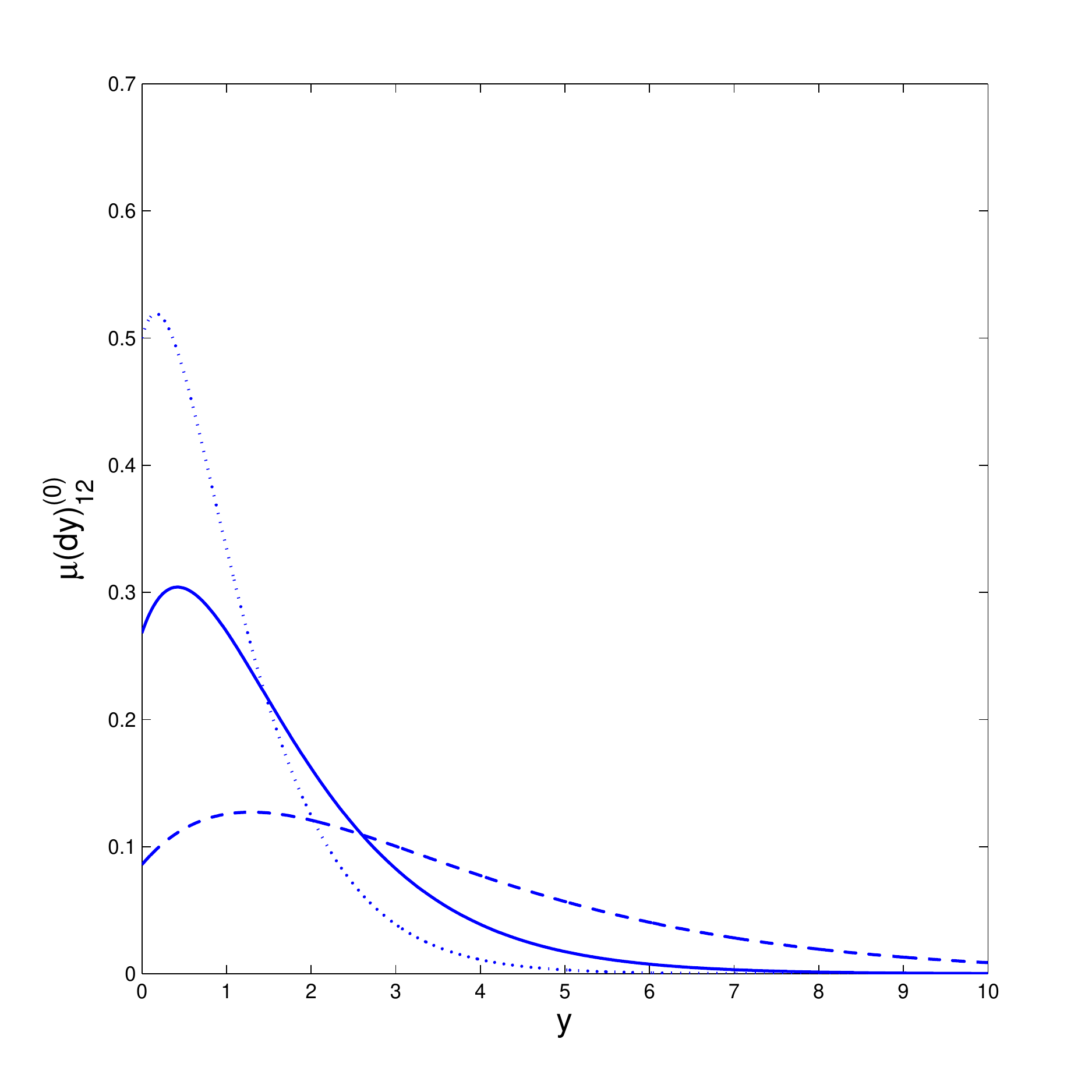}
\caption{The values of $\bmu(dy)^{(0)}_{11}/dy$ and $\bmu(dy)^{(0)}_{12}/dy$ in Example~\ref{ex1} for $b=1$, $a=4,\ 3,\ 2$ (dotted, solid, dashed line, respectively).}
\label{figex1}
\end{figure}	


\end{example}

We will now find Yaglom limit for strictly positive initial position of $X(0)=x>0$. Define matrices, for $n\geq 1$,
	\begin{eqnarray}
	\mathbf{W}(s^*, x-z)&=&\sum_{n=1}^\infty \frac{(x-z)^n}{n!}\sum_{i=1}^{n-1}\mathbf{D}(s^*)^i
	\times\mathbf{Q}_{21}(s^*)
	\mathbf{B}(s^*)
	\times\mathbf{D}(s^*)^{n-1-i},
	\nonumber\\
	\mathbf{W}_x(s^*)&=&\int_{z=0}^x \mathbf{W}(s^*, x-z)dz,
	\nonumber\\
	\mathbf{Z}_x(s^*,y)&=&\int_{z=0}^{\min\{x,y\}}\left(\mathbf{W}(s^*, x-z){\bf Q}_{21}(s^*)e^{{\bf K}(s^*)(y-z)}\right.
	\nonumber\\
	&&\left.+e^{{\bf D}(s^*)(x-z)}{\bf Q}_{21}(s^*)\mathbf{H}(s^*, y-z)
	\right)dz,
	\nonumber\\
	{\mathbf{Z}}_x(s^*)&=&\int_{y=0}^\infty \mathbf{Z}_x(s^*,y)\;dy,
	\end{eqnarray}
	and column vectors
	\begin{eqnarray}
		\widetilde{\mathbf{Z}}_x(s^*)
		&=&\widetilde{\mathbf{Z}}_x(s^*)_{11}{\bf 1}
		+ \widetilde{\mathbf{Z}}_x(s^*)_{12}{\bf 1}
		\nonumber\\
		&&
		+
		\left[
		\begin{array}{cc}
		\widetilde{\mathbf{Z}}_x(s^*)_{11} & 
		\widetilde{\mathbf{Z}}_x(s^*)_{12}
		\end{array}
		\right]
		\left[
		\begin{array}{c}
		{\bf T}_{10}\\
		{\bf T}_{20}
		\end{array}
		\right]
		(-({{\bf T}}_{00}-s^*{\bf I})^{-1})
		{\bf 1},
		\nonumber\\
		\widetilde{\widetilde{\mathbf{Z}}}_x(s^*)
		&=&\widetilde{\widetilde{\mathbf{Z}}}_x(s^*)_{21}{\bf 1}
		+ \widetilde{\widetilde{\mathbf{Z}}}_x(s^*)_{22}{\bf 1}
		\nonumber\\
		&&
		+
		\left[
		\begin{array}{cc}
		\widetilde{\widetilde{\mathbf{Z}}}_x(s^*)_{21} & 
		\widetilde{\widetilde{\mathbf{Z}}}_x(s^*)_x(s^*)_{22}
		\end{array}
		\right]
		\left[
		\begin{array}{c}
		{\bf T}_{10}\\
		{\bf T}_{20}
		\end{array}
		\right]
		(-({{\bf T}}_{00}-s^*{\bf I})^{-1})
		{\bf 1},
	\end{eqnarray}
	where
	\begin{eqnarray}
	\widetilde{\mathbf{Z}}_x(s^*)_{11}&=&
	{\mathbf{Z}}_x(s^*){\bf C}_1^{-1}
	\nonumber\\
	\widetilde{\mathbf{Z}}_x(s^*)_{12}&=&
	\left(
	 {\bf E}^{(x)}(s^*)_{21}{\bf C}_1
	 \mathbf{B}(s^*)
	 +{\mathbf{Z}}_x(s^*){\bf \Psi}(s^*)+\mathbf{W}_x(s^*)
	 \right) {\bf C}_2^{-1}
	 \nonumber\\
	 \widetilde{\widetilde{\mathbf{Z}}}_x(s^*)_{21}
	 &=&
	 \mathbf{B}(s^*){\bf E}^{(x)}(s^*)_{21}
	 +\mathbf{\Psi}(s^*)\mathbf{Z}_x(s^*){\bf C}_1^{-1}
	 +\mathbf{H}(s^*){\bf C}_1^{-1},
	 \nonumber\\
	 \widetilde{\widetilde{\mathbf{Z}}}_x(s^*)_{22}&=&
	 \mathbf{B}(s^*){\bf E}^{(x)}(s^*)_{22}
	 +\mathbf{\Psi}(s^*)
	 \Big( 
	 {\bf E}^{(x)}(s^*)_{21}{\bf C}_1\mathbf{B}(s^*){\bf C}_2^{-1}
	 +\mathbf{Z}_x(s^*){\bf\Psi}(s^*){\bf C}_2^{-1}
	 + \mathbf{W}(s^*){\bf C}_2^{-1}
	 \Big)
	 \nonumber\\
	 &&
	 +
	 \left(
	 (-{{\bf K}(s^*))^{-1}}\mathbf{B}(s^*)
	 +
	 {\bf H}(s^*){\bf\Psi}(s^*)
	 \right){\bf C}_2^{-1},
	 \nonumber\\
	\end{eqnarray}
with ${\bf E}^{(x)}(s^*)_{21}=\int_{y=0}^{\infty}{\bf E}(dy)^{(x)}(s^*)_{21}$ as considered in Remark~\ref{howtogetEx21s}, and
\begin{eqnarray}
{\bf E}^{(x)}(s^*)_{22}&=&
{\bf E}^{(x)}(s^*)_{21}{\bf C}_1{\bf\Psi}(s^*){\bf C}_2^{-1}
+\int_{y=0}^xe^{{\bf D}(s^*)(x-y)}{\bf C}_2^{-1}dy
\nonumber\\
&=&
{\bf E}^{(x)}(s^*)_{21}{\bf C}_1{\bf\Psi}(s^*){\bf C}_2^{-1}
-\int_{w=0}^xe^{{\bf D}(s^*)w}{\bf C}_2^{-1}dw
\nonumber\\
&=&
{\bf E}^{(x)}(s^*)_{21}{\bf C}_1{\bf\Psi}(s^*){\bf C}_2^{-1}
-
({\bf D}(s^*))^{-1}
\left(
e^{{\bf D}(s^*)x}-{\bf I}
\right){\bf C}_2^{-1}
.
\end{eqnarray}

\begin{theorem}\label{thm:main2}
For $x>0$ the matrix $\bmu(dy)^{(x)}$ is unique and
\begin{eqnarray*}
\bmu(dy)^{(x)}_{21}&=& 
diag({\widetilde{\bf Z}_x(s^*)})^{-1}
{{\bf Z}_x(s^*,y){\bf C}_1^{-1}}
dy,\nonumber\\
\bmu(dy)^{(x)}_{22}&=&
diag({\widetilde{\bf Z}_x(s^*)})^{-1}
\Big( 
{\bf E}(dy)^{(x)}(s^*)_{21}{\bf C}_1\mathbf{B}(s^*){\bf C}_2^{-1}
+\mathbf{Z}_x(s^*,y){\bf\Psi}(s^*){\bf C}_2^{-1}dy 
\nonumber\\
&&
+ \mathbf{W}(s^*, x-y){\bf 1}\{y<x\}{\bf C}_2^{-1}dy
\Big)
,\label{qs21x}\\
\bmu(dy)^{(x)}_{20}&=&
\left[
\begin{array}{cc}
\bmu(dy)^{(x)}_{21}&\bmu(dy)^{(x)}_{22}
\end{array}
\right]
\left[
\begin{array}{c}
{\bf T}_{10}\\
{\bf T}_{20}
\end{array}
\right]
(-({{\bf T}}_{00}-s^*{\bf I})^{-1}),\nonumber
\\
\bmu(dy)^{(x)}_{11}&=&
diag({\widetilde{\widetilde{\bf Z}}_x(s^*)})^{-1}
\Big(
\mathbf{B}(s^*){\bf E}(dy)^{(x)}(s^*)_{21}
+\mathbf{\Psi}(s^*)\mathbf{Z}_x(s^*,y){\bf C}_1^{-1} dy
\nonumber\\
&&
+\mathbf{H}(s^*, y-x){\bf C}_1^{-1}{\bf 1}\{y>x\}dy
\Big)
,
\nonumber\\
\bmu(dy)^{(x)}_{12}&=&
diag({\widetilde{\widetilde{\bf Z}}_x(s^*)})^{-1}
\Big\{
\Big(
\mathbf{B}(s^*)
{\bf E}(dy)^{(x)}(s^*)_{22}
+\mathbf{\Psi}(s^*)
\Big( 
{\bf E}(dy)^{(x)}(s^*)_{21}{\bf C}_1
\mathbf{B}(s^*)
{\bf C}_2^{-1}
\nonumber\\
&&
+\mathbf{Z}_x(s^*,y){\bf\Psi}(s^*){\bf C}_2^{-1}dy
+ \mathbf{W}(s^*, x-y){\bf 1}\{y<x\}{\bf C}_2^{-1}dy
\Big)
\Big)
\nonumber\\
&&
+
\left(
e^{{\bf K}(s^*)(y-x)}\mathbf{B}(s^*)
+
{\bf H}(s^*,y-x){\bf\Psi}(s^*)
\right){\bf C}_2^{-1}{\bf 1}\{y>x\}dy
\Big\},
\nonumber\\
\bmu(dy)^{(x)}_{10}&=&
\left[
\begin{array}{cc}
\bmu(dy)^{(x)}_{11}&\bmu(dy)^{(x)}_{12}
\end{array}
\right]
\left[
\begin{array}{c}
{\bf T}_{10}\\
{\bf T}_{20}
\end{array}
\right]
(-({{\bf T}}_{00}-s^*{\bf I})^{-1}).
\label{eq:Edmodyfikacja}
\end{eqnarray*}

\end{theorem}

From Theorems \ref{thm:main} and \ref{thm:main2} it follows the following corollary.
\begin{corollary}
Yaglom limit depends on the initial position of the fluid level $X(0)=x$ in the model.
\end{corollary}

\begin{remark}
There has been a conjecture that
Yaglom limit does not depend on initial position of the Markov process. However, a counterexample to this
conjecture was already demonstrated by Foley and McDonald \cite{FolMcDon}.
Our model produces another example of the same kind.
\end{remark}

{\bf Proof:}
Our proof is again based on Theorem~\ref{t.tauberian} and \eqref{mu_fraction}.
Note that
\begin{eqnarray}
\lefteqn{e^{\mathbf{D}(s)(x-z)}=\lim_{K\to+\infty}\sum_{n=0}^K
\frac{(x-z)^n}{n!}\left(\mathbf{Q}_{22}(s)+\mathbf{Q}_{21}(s)\mathbf{\Psi}(s)\right)^n}\nonumber\\
&&\quad=\lim_{K\to+\infty}\sum_{n=0}^k\frac{(x-z)^n}{n!}\left(\mathbf{Q}_{22}(s^*)+\mathbf{Q}_{21}(s^*)
(\mathbf{\Psi}(s^*)
-\mathbf{B}(s^*)
\sqrt{s-s^*}+o(\sqrt{s-s^*})))\right)^n\nonumber\\
&&\quad= e^{\mathbf{D}(s^*)(x-z)}
- \sqrt{s-s^*}
\mathbf{W}(s^*, x-z)+o(\sqrt{s-s^*}).\nonumber\\
&&\label{edoD}
\end{eqnarray}
By \eqref{edoK}, \eqref{edoD}, Lemmas~\ref{QS_form} and \ref{th:Edy} and Theorem~\ref{lem:psi_form}, we have
\begin{eqnarray}
{\bf E}(dy)^{(x)}(s)_{21}&=&
\int_{z=0}^{\min\{x,y\}}
e^{{\bf D}(s)(x-z)}{\bf Q}_{21}(s)e^{{\bf K}(s)(y-z)}{\bf C}_1^{-1}dzdy
\nonumber\\
&=&
\int_{z=0}^{\min\{x,y\}}\left(e^{\mathbf{D}(s^*)(x-z)}
- \sqrt{s-s^*}
\mathbf{W}^*(s^*, x-z)\right)\mathbf{Q}_{21}(s^*)
\nonumber\\
&&\times \left(e^{\mathbf{K}(s^*)(y-z)}
- \sqrt{s-s^*}
\mathbf{H}(s^*, y-z)\right){\bf C}_1^{-1}dz dy
+o(\sqrt{s-s^*})
\nonumber\\
&=&
{\bf E}(dy)^{(x)}(s^*)_{21}
- \sqrt{s-s^*}
\mathbf{Z}_x(s^*,y){\bf C}_1^{-1} dy +o(\sqrt{s-s^*}),
\end{eqnarray}
and
\begin{eqnarray}
\lefteqn{
{\bf E}(dy)^{(x)}(s)_{22}=
{\bf E}(dy)^{(x)}(s)_{21}{\bf C}_1{\bf\Psi}(s){\bf C}_2^{-1}
+e^{{\bf D}(s)(x-y)}{\bf C}_2^{-1}{\bf 1}\{y<x\}dy
}
\nonumber\\
&=&
\left(
{\bf E}(dy)^{(x)}(s^*)_{21}{\bf C}_1
- \sqrt{s-s^*}
\mathbf{Z}_x(s^*,y) dy
\right)
\left(
{\bf\Psi}(s^*)-{\bf B}(s^*)\sqrt{s-s^*}
\right)
{\bf C}_2^{-1}
\nonumber\\
&&
+
\left(
e^{\mathbf{D}(s^*)(x-y)}
- \sqrt{s-s^*}
\mathbf{W}(s^*, x-y)
\right){\bf C}_2^{-1}{\bf 1}\{y<x\}dy+o(\sqrt{s-s^*})
\nonumber\\
&=&{\bf E}(dy)^{(x)}(s^*)_{22}
- \sqrt{s-s^*}
\Big( 
{\bf E}(dy)^{(x)}(s^*)_{21}{\bf C}_1
(-\mathbf{B}(s^*))
{\bf C}_2^{-1}
+\mathbf{Z}_x(s^*,y){\bf\Psi}(s^*){\bf C}_2^{-1}dy
\nonumber\\
&&
+ \mathbf{W}(s^*, x-y){\bf 1}\{y<x\}{\bf C}_2^{-1}dy
\Big)
+o(\sqrt{s-s^*}),
\nonumber\\
\end{eqnarray}
and
\begin{eqnarray}
\lefteqn{
{\bf E}(dy)^{(x)}(s)_{20}
=
\left[
\begin{array}{cc}
{\bf E}(dy)^{(x)}(s)_{21}&{\bf E}(dy)^{(x)}(s)_{22}
\end{array}
\right]
\left[
\begin{array}{c}
{\bf T}_{10}\\
{\bf T}_{20}
\end{array}
\right]
(-({{\bf T}}_{00}-s{\bf I})^{-1})}
\nonumber\\
&=&
{\bf E}(dy)^{(x)}(s^*)_{20}
\nonumber\\
&&
- \sqrt{s-s^*}
	\left[
	\begin{array}{cc}
	{\mathbf{Z}}_x(s^*,y) & 
	{\bf E}(dy)^{(x)}(s^*)_{21}{\bf C}_1
	\mathbf{B}(s^*)
	+{\mathbf{Z}}_x(s^*){\bf \Psi}(s^*)dy+\mathbf{W}(s^*,x-y){\bf 1}\{y<x\}dy
	\end{array}
	\right]
			\nonumber\\
			&&
	\times
	\left[
	\begin{array}{c}
	{\bf C}_1^{-1}{\bf T}_{10}\\
	{\bf C}_2^{-1}{\bf T}_{20}
	\end{array}
	\right]
	(-({{\bf T}}_{00}-s^*{\bf I})^{-1})
+o(\sqrt{s-s^*}),
\nonumber\\
\end{eqnarray}
and
\begin{eqnarray}
{\bf E}^{(x)}(s)_2&=&
\int_{y=0}^{\infty}{\bf E}(dy)^{(x)}(s)_{21}{\bf 1}
+\int_{y=0}^{\infty}{\bf E}(dy)^{(x)}(s)_{22}{\bf 1}
+\int_{y=0}^{\infty}{\bf E}(dy)^{(x)}(s)_{20}{\bf 1}
\nonumber\\
&=&
{\bf E}^{(x)}(s^*)_1
- \sqrt{s-s^*}\widetilde{\bf Z}_x(s^*)
+o(\sqrt{s-s^*}).
\nonumber\\
\end{eqnarray}
Thus the expressions for $\bmu(dy)^{(x)}_{21}$, $\bmu(dy)^{(x)}_{22}$ and $\bmu(dy)^{(x)}_{20}$ follow by argument similar to the proof of Theorem~\ref{thm:main}.

Furthermore, by \eqref{edoK}, Lemmas~\ref{QS_form} and \ref{th:Edy} and Theorem~\ref{lem:psi_form}, we have
\begin{eqnarray}
\lefteqn{
{\bf E}(dy)^{(x)}(s)_{11}={\bf\Psi}(s){\bf E}(dy)^{(x)}(s)_{21}
+e^{{\bf K}(s)(y-x)}{\bf C}_1^{-1}{\bf 1}\{y>x\}dy
}
\nonumber\\
&=&
(\mathbf{\Psi}(s^*)
-\mathbf{B}(s^*)
\sqrt{s-s^*})
({\bf E}(dy)^{(x)}(s^*)_{21}
- \sqrt{s-s^*}
\mathbf{Z}_x(s^*,y){\bf C}_1^{-1} dy)
\nonumber\\
&&
+\left(e^{\mathbf{K}(s^*)(y-x)}
- \sqrt{s-s^*}
\mathbf{H}(s^*, y-x)\right)
{\bf C}_1^{-1}{\bf 1}\{y>x\}dy
+o(\sqrt{s-s^*})
\nonumber\\
&=&{\bf E}(dy)^{(x)}(s^*)_{11}
- \sqrt{s-s^*}
\Big(
\mathbf{B}(s^*){\bf E}(dy)^{(x)}(s^*)_{21}
+\mathbf{\Psi}(s^*)\mathbf{Z}_x(s^*,y){\bf C}_1^{-1} dy
\nonumber\\
&&
+\mathbf{H}(s^*, y-x){\bf C}_1^{-1}{\bf 1}\{y>x\}dy
\Big)
+o(\sqrt{s-s^*}),
\nonumber\\
\end{eqnarray}
and
\begin{eqnarray}
\lefteqn{
{\bf E}(dy)^{(x)}(s)_{12}
=
{\bf\Psi}(s){\bf E}(dy)^{(x)}(s)_{22}
+e^{{\bf K}(s)(y-x)}{\bf\Psi}(s){\bf C}_2^{-1}{\bf 1}\{y> x\}dy
}
\nonumber\\
&=&
\left(\mathbf{\Psi}(s^*)
-\mathbf{B}(s^*)\sqrt{s-s^*}
\right)
\Big( {\bf E}(dy)^{(x)}(s^*)_{22}
- \sqrt{s-s^*}
\Big( 
{\bf E}(dy)^{(x)}(s^*)_{21}{\bf C}_1
\mathbf{B}(s^*)
{\bf C}_2^{-1}
\nonumber\\
&&
+\mathbf{Z}_x(s^*,y){\bf\Psi}(s^*){\bf C}_2^{-1}dy
+ \mathbf{W}(s^*, x-y){\bf 1}\{y<x\}{\bf C}_2^{-1}dy
\Big)
\nonumber\\
&&
+
\Big(
e^{{\bf K}(s^*)(y-x)}{\bf\Psi}(s^*){\bf C}_2^{-1}{\bf 1}\{y> x\}dy
\nonumber\\
&&
- \sqrt{s-s^*}
\left(
e^{{\bf K}(s^*)(y-x)}\mathbf{B}(s^*)
+
{\bf H}(s^*,y-x){\bf\Psi}(s^*)
\right){\bf C}_2^{-1}dy
\Big)
+o(\sqrt{s-s^*})
\nonumber\\
&=&
{\bf E}(dy)^{(x)}(s^*)_{12}
- \sqrt{s-s^*}
\Big(
\mathbf{B}(s^*)
{\bf E}(dy)^{(x)}(s^*)_{22}
+\mathbf{\Psi}(s^*)
\Big( 
{\bf E}(dy)^{(x)}(s^*)_{21}{\bf C}_1
\mathbf{B}(s^*)
{\bf C}_2^{-1}
\nonumber\\
&&
+\mathbf{Z}_x(s^*,y){\bf\Psi}(s^*){\bf C}_2^{-1}dy
+ \mathbf{W}(s^*, x-y){\bf 1}\{y<x\}{\bf C}_2^{-1}dy
\Big)
\Big)
\nonumber\\
&&
- \sqrt{s-s^*}
\left(
e^{{\bf K}(s^*)(y-x)}\mathbf{B}(s^*)
+
{\bf H}(s^*,y-x){\bf\Psi}(s^*)
\right){\bf C}_2^{-1}{\bf 1}\{y>x\}dy
+o(\sqrt{s-s^*})
.
\nonumber\\
\end{eqnarray}
Thus the expressions for $\bmu(dy)^{(x)}_{11}$, $\bmu(dy)^{(x)}_{12}$ and $\bmu(dy)^{(x)}_{10}$ follow by a similar argument, with
\begin{eqnarray}
{\bf E}^{(x)}(s)_1&=&
\int_{y=0}^{\infty}{\bf E}(dy)^{(x)}(s)_{11}{\bf 1}
+\int_{y=0}^{\infty}{\bf E}(dy)^{(x)}(s)_{12}{\bf 1}
+\int_{y=0}^{\infty}{\bf E}(dy)^{(x)}(s)_{10}{\bf 1}
\nonumber\\
&=&
{\bf E}^{(x)}(s^*)_1
- \sqrt{s-s^*}\widetilde{\widetilde{\bf Z}}_x(s^*)
+o(\sqrt{s-s^*}).
\nonumber\\
\end{eqnarray}
 \rule{9pt}{9pt}

\section{Example with non-scalar ${\bf\Psi}(s)$}
\label{sec:s_ex}

Below we construct an example where, unlike in Example~\ref{ex1}, key quantities are matrices, rather than scalars. We derive expressions for this example analytically and illustrate these results with some numerical output as well.

\begin{example}\label{ex2} \rm Consider a system with $N=2$ sources based on example analysed in~\cite{AMS}. Let $\mathcal{S}=\{1,2,3\}$, $\mathcal{S}_1=\{1\}$, $\mathcal{S}_2=\{2,3\}$, $c_1=1$, $c_2=c_3=-1$, and
\begin{eqnarray}
{\bf T}
&=&\left[
\begin{array}{cc}
{\bf T}_{11}&{\bf T}_{12}\\
{\bf T}_{21}&{\bf T}_{22}
\end{array}
\right]
=\left[
\begin{array}{c|cc}
-2\lambda&2\lambda&0\\
\hline
1&-(1+\lambda)&\lambda\\
0&2&-2
\end{array}
\right],
\nonumber\\
{\bf Q}(s)
&=&\left[
\begin{array}{cc}
{\bf Q}_{11}(s)&{\bf Q}_{12}(s)\\
{\bf Q}_{21}(s)&{\bf Q}_{22}(s)
\end{array}
\right]
=\left[
\begin{array}{c|cc}
-(2\lambda +s)&2\lambda&0\\
\hline
1&-(1+\lambda +s)&\lambda\\
0&2&-(2+s)
\end{array}
\right],
\nonumber
\end{eqnarray}
with some parameter $\lambda>\sqrt{2}-1$ so that the process is stable. In our plots of the output below, we will assume the value $\lambda=2.5$.

Denote by $[x\ z]={\bf\Psi}(s)=\int_{t=0}^{\infty}e^{-st}\bpsi(t)dt$ the minimum nonnegative solution of~(\ref{Ric}), here equivalent to
\begin{eqnarray}
[0\ 0]&=&[2\lambda\ 0]-(2\lambda+s)[x\ z]
+[x\ z]
\left[
\begin{array}{cc}
-(1+\lambda+s)&\lambda\\
2&-(2+s)
\end{array}
\right]
\nonumber\\
&&
+[x\ z]
\left[
\begin{array}{c}
1\\
0
\end{array}
\right]
[x\ z],
\end{eqnarray}
which we write as a system of equations
\begin{eqnarray}
0&=&x^2-(1+3\lambda+2s)x+2z+2\lambda,\label{eq:scalar1}\\
0&=&-(2+2\lambda+2s-x)z+\lambda x.
\label{eq:scalar}
\end{eqnarray}

The minimum nonnegative solution $[x\ z]$ of~\eqref{eq:scalar1}-\eqref{eq:scalar} must be strictly positive, satisfy $2+2\lambda+2s-x> 0$, and occur at the intersection of the two curves,
\begin{eqnarray}
z=z_1(x,s)&=&-\frac{1}{2}x^2+\frac{1}{2}(1+3\lambda+2s)x-\lambda,\label{z1}\\
z=z_2(x,s)&=&\lambda x/(2+2\lambda+2s-x).\label{z}
\end{eqnarray}

We consider the shape of the curves in~\eqref{z1}-\eqref{z} to facilitate the analysis that follows, see Figure~\ref{twocruves1}. It is a straightforward exercise to verify that, when $s=0$, we have $z_1(x,0)<-\lambda<z_2(x,0)$ for all $x< 0$, and so the two curves may only intersect at some point $(x,z)$ with $x>0$.
\begin{figure}[htbp]
\centering
 \begin{tabular}{@{}cc@{}}
    \includegraphics[width=.4\textwidth]{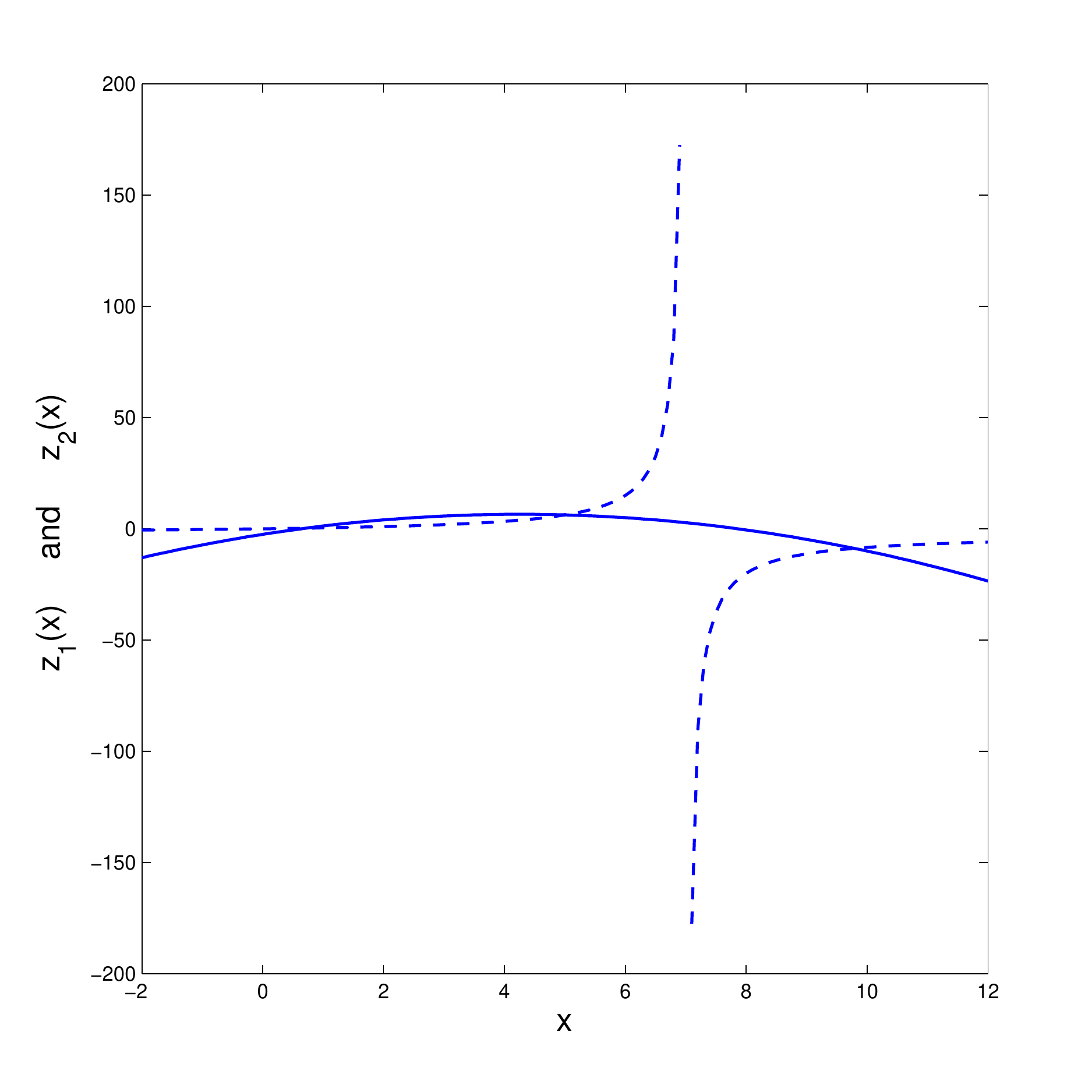} &
    \includegraphics[width=.4\textwidth]{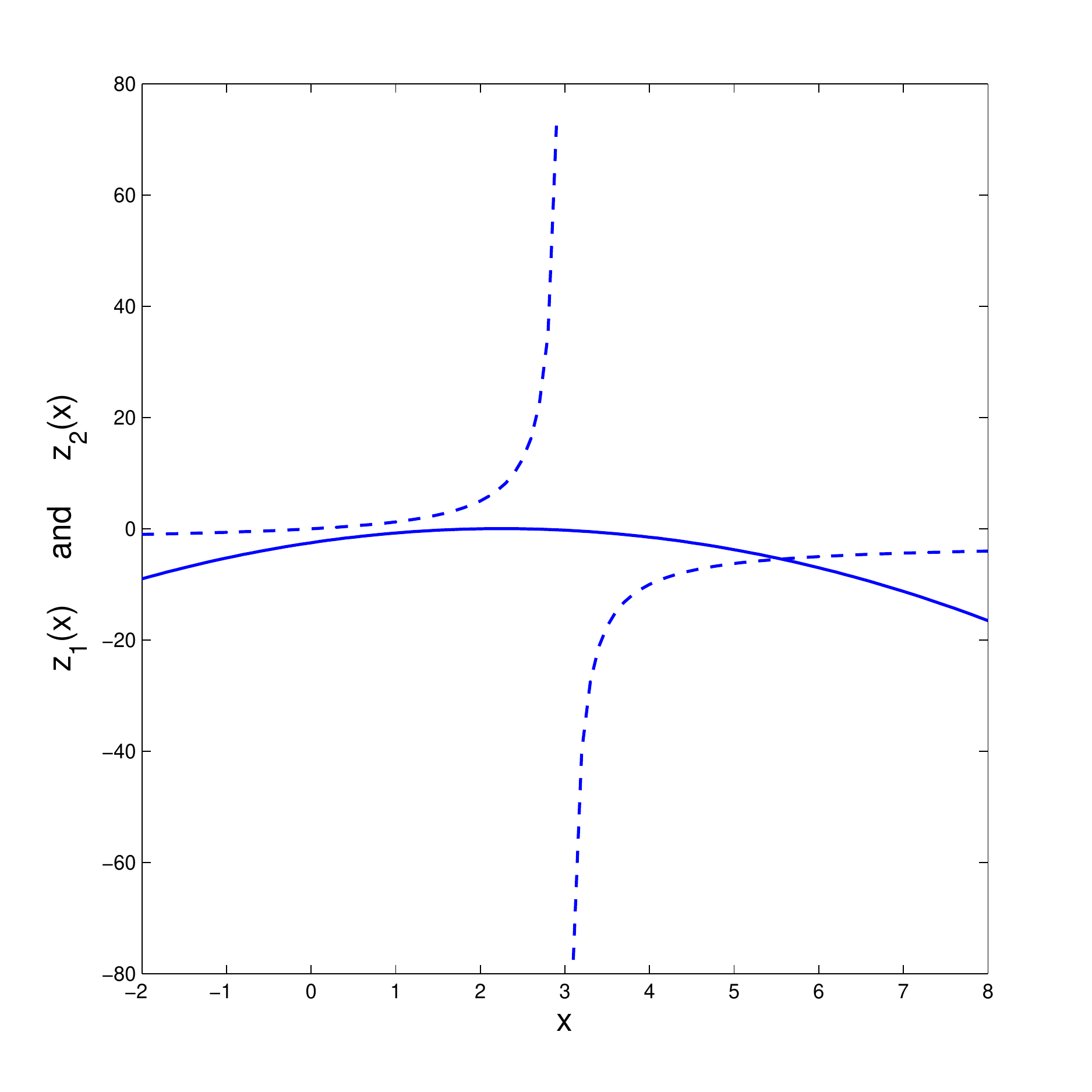} \\
  \end{tabular}
\caption{The plot of~\eqref{z1}-\eqref{z} for $s=0$ (left) and $s=-2$ (right), when $\lambda=2.5$.
}
\label{twocruves1}
\end{figure}

Further, when $2+2\lambda+2s-x> 0$, we have
\begin{eqnarray}
\frac{\partial z_2(x,s)}{\partial x}=
\frac{\lambda(2+2\lambda+2s-x)+\lambda x}
{(2+2\lambda+2s-x)^2}
&>&0,
\end{eqnarray}
and so,  when $s=0$, then the minimum nonnegative solution $[x\ z]$ of~\eqref{eq:scalar1}-\eqref{eq:scalar} is in fact the minimum real-valued solution of~\eqref{eq:scalar1}-\eqref{eq:scalar}.

Also, when $x>0$ and $2+2\lambda+2s-x> 0$, we have
\begin{eqnarray}
\frac{\partial z_1(x,s)}{\partial s}=x
&>&0
,\nonumber\\
\frac{\partial z_2(x,s)}{\partial s}=\frac{-2\lambda x}{(2+2\lambda+2s-x)^2}
&<&0,
\end{eqnarray}
and so as $s\downarrow s^*$ we have $z_1(x,s)\downarrow$ while $z_2(x,s)\uparrow$, until the two curves touch when $s=s^*$, and then move apart when $s<s^*$. Therefore, by the continuity of ${\bf\Psi}(s)$ argument as used in the proof of Lemma~\ref{lem:sstardeltastar2}, for all $s\in [s^*,0]$, ${\bf\Psi}(s)=[x\ z]$ is the minimum real-valued solution of~\eqref{eq:scalar1}-\eqref{eq:scalar}.

Instead of looking at the problem as two intersecting curves $z_1(x,s)$ and $z_2(x,s)$, we now look at it as one cubic curve $g_s(x)$. Substitute~\eqref{z} into~\eqref{eq:scalar1} and multiply by $(2+2\lambda+2s-x)$, to get
\begin{eqnarray}
0&=&
-x^3
+(3+5\lambda+4s)x^2
-(2+2\lambda+2s)(1+3\lambda+2s)x+(2+2\lambda+2s)2\lambda
\nonumber\\
&=&g_s(x),
\label{eq:cubic}
\end{eqnarray}
which is of the form
\begin{equation}\label{eq:cubicpar}
ax^3+bx^2+cx+d=0,
\end{equation}
with $g_s(0)=d>0$ (we have $d>0$ since $0<x<2+2\lambda+2s$ due to $z>0$ in~\eqref{z}). See the plots of $g_s(x)$ in Figure~\ref{plotg} for the case $\lambda=2.5$. Noting that $a=-1<0$, we conclude that when $s=s^*$, the the solution $[x\ z]$ corresponds to the local minimum,
\begin{equation}\label{eq:x}
x=\min \left\{ \frac{-b+ \sqrt{b^2-3ac}}{3a},\frac{-b - \sqrt{b^2-3ac}}{3a} \right\}
=\frac{-b+ \sqrt{b^2-3ac}}{3a},
\end{equation}
where
\begin{equation}\label{eq:crit}
b^2-3ac>0.
\end{equation}

\begin{figure}[htbp]
	\centering
	\begin{tabular}{@{}cc@{}}
		\includegraphics[width=.4\textwidth]{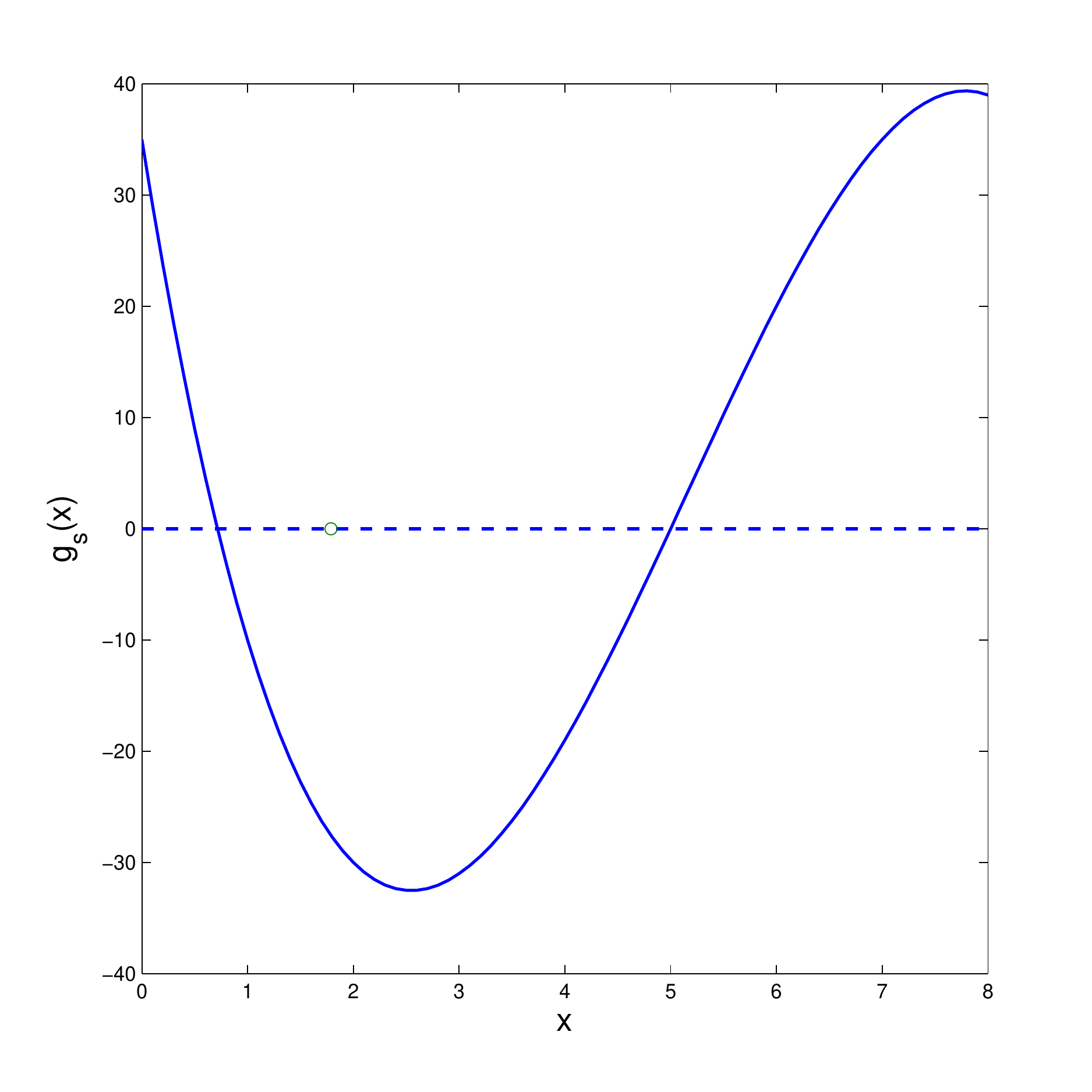} &
		\includegraphics[width=.4\textwidth]{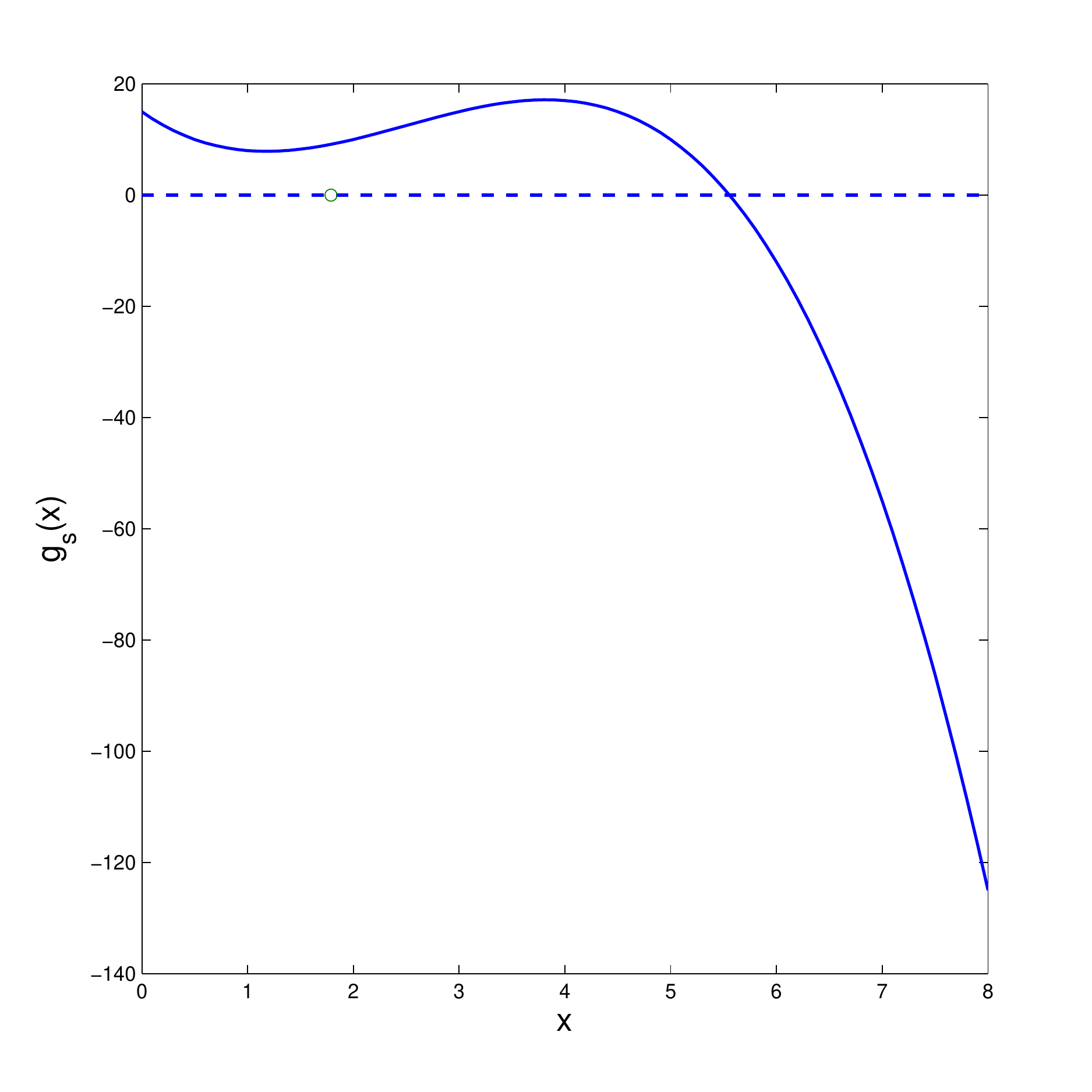} \\
	\end{tabular}
	\includegraphics[scale=0.4]{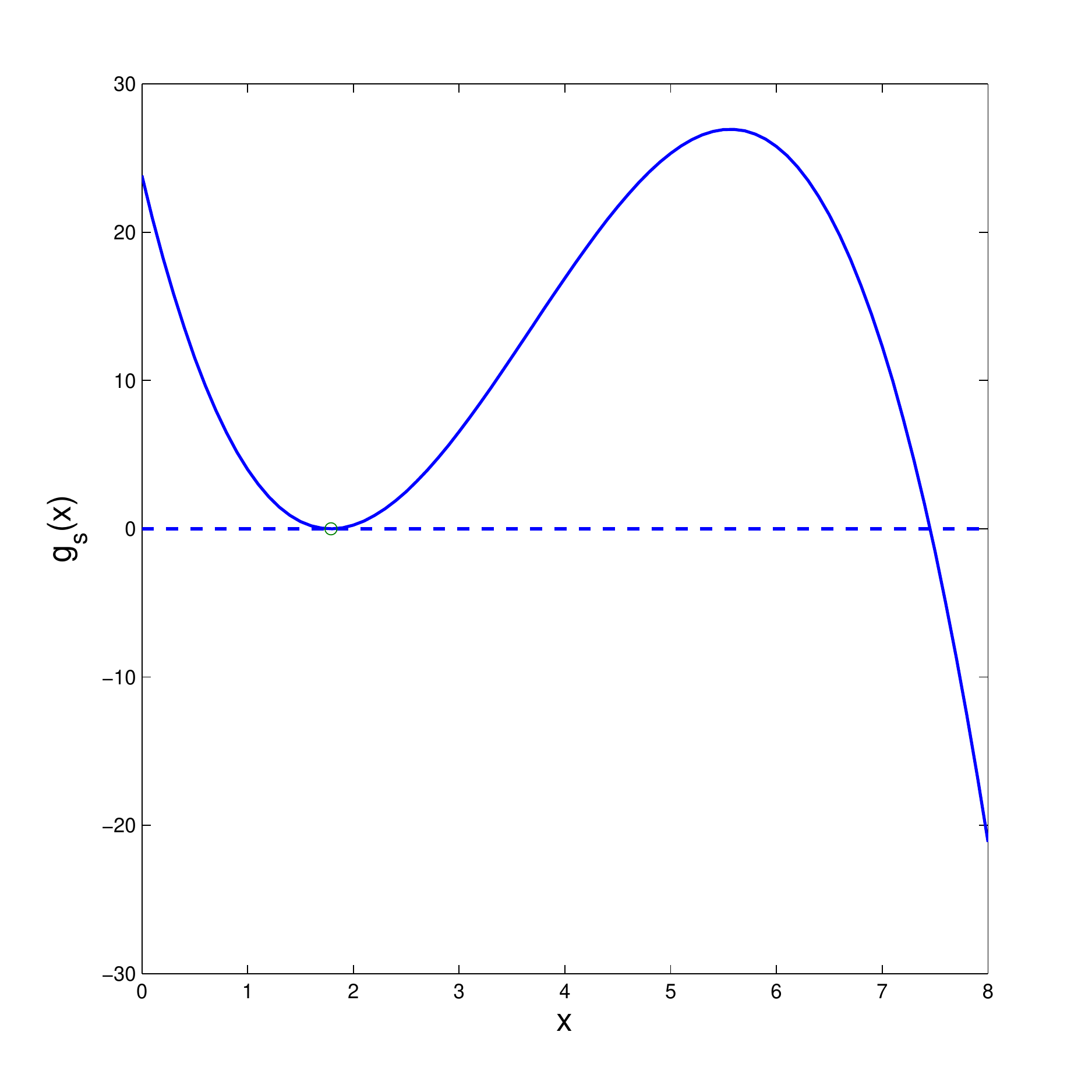}
	\caption{The plot of~\eqref{eq:cubic} for $s=0$ (top left) and $s=-2$ (top right) and $s=-1.1178$, when $\lambda=2.5$.
	}
	\label{plotg}
\end{figure}

We transform the cubic equation~\eqref{eq:cubicpar} into
\begin{equation}\label{eq_y}
y^3+py+q=0
\end{equation}
using
\begin{eqnarray}
x&=&y-\frac{b}{3a}
,\label{xy}
\end{eqnarray}
with
\begin{eqnarray}
p&=&\frac{3ac-b^2}{3a^2}
= s\times c_p^{(1)}+s^2\times c_p^{(2)}+c_p
,\label{p}
\end{eqnarray}
for suitable $c_p^{(1)}$, $c_p^{(2)}$, $c_p$, and
\begin{eqnarray}
q&=&
\frac{2b^3+27a^2d-9abc}{27a^3}
= s\times c_q^{(1)}+s^2\times c_q^{(2)}+s^3\times c_q^{(3)}
+c_q
,\label{q}
\end{eqnarray}
for suitable $c_q^{(1)}$, $c_q^{(2)}$, $c_q^{(3)}$, $c_q$.

Below, we choose the convention that we write $p(s)$ to demonstrate the $p$ is a function of $s$, with similar notation applied for other quantities like $q$, $x$, $y$ and so on. Observe that
\begin{eqnarray*}
s^3-(s^*)^3&=&(s-s^*)(s^2+ss^*+(s^*)^2)
=C_3\times (s-s^*)+o(s-s^*),\\
s^2-(s^*)^2&=&(s-s^*)(s+s^*)
=C_2\times (s-s^*)+o(s-s^*),
\end{eqnarray*}
where $C_3=3(s^*)^2$ and $C_2=2s^*$, and so by~\eqref{p}-\eqref{q},
\begin{eqnarray}
p(s)-p(s^*)&=&C_p\times (s-s^*)+o(s-s^*),\label{trreprp}\\
q(s)-q(s^*)&=&C_q\times (s-s^*)+o(s-s^*),\label{trreprq}
\end{eqnarray}
where constants $C_p$ and $C_q$ are given by
\begin{eqnarray}
C_p&=& c_p^{(1)}+C_2\times c_p^{(2)}
,\nonumber\\
C_q&=&c_q^{(1)}+C_2\times c_q^{(2)}+C_3\times c_q^{(3)}
.
\end{eqnarray}

Consider~\eqref{eq_y} and apply Vi\'eta's substitution,
\begin{equation}\label{Viet}
y=u-\frac{p}{3u},
\end{equation}
where $u^3$ solves the quadratic equation,
\begin{equation}\label{eq:u}
(u^3)^2+qu^3-\frac{p^3}{27}=0,
\end{equation}
and the two solutions are
\begin{equation}
u^3(s)=\frac{-q(s)\pm \sqrt{\Delta (s)}}{2},
\end{equation}
with
\begin{equation}\label{eq:Deltas}
\Delta (s)= q^2(s)+4\times \frac{p^3(s)}{27},
\end{equation}
where $\Delta(s)<0$ for $s>s^*$ and the repeated root requires
\begin{equation}\label{sstarex2}
\Delta(s^*)=q^2(s^*)+4\times \frac{p^3(s^*)}{27}=0.
\end{equation}

When $s>s^*$, the three (real) solutions of~\eqref{eq_y} are the three cubic roots,
\begin{eqnarray}
y_0,y_1,y_2&=&
\left(
\frac{-q(s) + \sqrt{\Delta (s)}}{2}
\right)^{1/3},
\end{eqnarray}
and we choose the minimum
\begin{equation}
y(s)=\min\{y_0(s),y_1(s),y_2(s)\},
\end{equation}
which corresponds to the minimum $x(s)=\mathbf{\Psi}(s)_1$ where $\mathbf{\Psi}(s)_i$ denotes $i$th element of $\mathbf{\Psi}(s)$.

Therefore, by~\eqref{trreprq},
\begin{eqnarray}
u^3(s)-u^3(s^*)&=&
\frac{-q(s) + \sqrt{\Delta (s)}}{2}
+ \frac{q(s^*)}{2}
\nonumber\\
&=&
-\frac{1}{2}C_q\times (s-s^*) + \frac{1}{2}\sqrt{\Delta(s)}
+o(s-s^*).
\end{eqnarray}

Now,
\begin{eqnarray}
\lefteqn{
\Delta(s)=\Delta(s)-\Delta(s^*)
}
\nonumber\\
&=&
\frac{1}{2}(q(s)-q(s^*))(q(s)+q(s^*))
+\frac{4}{27} (p(s)-p(s^*))(p(s)^2+p(s)p(s^*)+p(s^*)^2),
\nonumber\\
\end{eqnarray}
and so by \eqref{trreprp}-\eqref{trreprq},
\begin{equation}\label{u3}
\Delta(s)= C_\Delta\times (s-s^*)+o(s-s^*)
\end{equation}
where the constant $C_\Delta<0$ is given by
\begin{eqnarray}
C_\Delta &=&
\frac{1}{2}C_q\times 2q(s^*)
+\frac{4}{27} C_p\times 3p^2(s^*).
\end{eqnarray}

Therefore,
\begin{eqnarray*}
\lim_{s\downarrow s^*}
\left(
\frac{u^3(s)-u^3(s^*)}{\sqrt{s-s^*}}
\right)
&=&
\lim_{s\downarrow s^*}
\left(
-\frac{1}{2} C_q\sqrt{s-s^*}
+
\frac{1}{2}\sqrt{\frac{C_\Delta\times (s-s^*)+o(s-s^*)}{s-s^*}}
\right)
\\
&=&
\frac{1}{2}\sqrt{C_\Delta},
\end{eqnarray*}
and
\begin{eqnarray*}
\lim_{s\downarrow s^*}
\left(
\frac{u(s)-u(s^*)}{\sqrt{s-s^*}}
\right)
&=&
\lim_{s\downarrow s^*}
\left(
\frac{u(s)-u(s^*)}{\sqrt{s-s^*}}
\times
\frac{u^2(s)+u(s)u(s^*)+u^2(s^*)}{u^2(s)+u(s)u(s^*)+u^2(s^*)}
\right)
\nonumber\\
&=&
\lim_{s\downarrow s^*}
\left(
\frac{u^3(s)-u^3(s^*)}{\sqrt{s-s^*}}
\times
\frac{1}{u^2(s)+u(s)u(s^*)+u^2(s^*)}
\right)
\nonumber\\
&=&\frac{1}{6u^2(s^*)}\sqrt{C_\Delta},
\end{eqnarray*}
where $u(s^*)\neq 0$ by~\eqref{eq:u}, since $p(s^*)\neq 0$ due to~\eqref{eq:crit} and~\eqref{p}.

From the above we conclude that by~\eqref{Viet},
\begin{eqnarray*}
	\lefteqn{
\lim_{s\downarrow s^*}
\left(
\frac{y(s)-y(s^*)}{\sqrt{s-s^*}}
\right)
=
\lim_{s\downarrow s^*}
\left(
\frac{u(s)-u(s^*)}{\sqrt{s-s^*}} -\frac{1}{3\sqrt{s-s^*}}
\left(\frac{p(s)}{u(s)}-
\frac{p(s^*)}
{u(s^*)}\right)
\right)
}\\
&=&
\lim_{s\downarrow s^*}
\left(
\frac{u(s)-u(s^*)}{\sqrt{s-s^*}}
-\frac{(p(s)-p(s^*))u(s^*)}{3\sqrt{s-s^*}u(s)u(s^*)}
+\frac{p(s^*)(u(s)-u(s^*))}{3\sqrt{s-s^*}u(s)u(s^*)}
\right)
\\
&=&
\pm \frac{1}{6u^2(s^*)}\sqrt{C_\Delta}
-0
 \frac{1}{6u^2(s^*)}\sqrt{C_\Delta}\frac{p(s^*)}{3u^2(s^*)}
\\
&=&  \frac{1}{6u^2(s^*)}\sqrt{C_\Delta}\left(1+\frac{p(s^*)}{3u^2(s^*)}\right). \end{eqnarray*}

Therefore, by~\eqref{xy}, we have
\begin{eqnarray}
\lim_{s\downarrow s^*}
\left(
\frac{\mathbf{\Psi}(s)_1-\mathbf{\Psi}(s^*)_1}{\sqrt{s-s^*}}
\right)
&=&
\lim_{s\downarrow s^*}
\left(
\frac{x(s)-x(s^*)}{\sqrt{s-s^*}}
\right)
\nonumber\\
&=&
\lim_{s\downarrow s^*}
\left(
\frac{y(s)-y(s^*)}{\sqrt{s-s^*}} +o(1)
\right)
\nonumber\\
&=&
\frac{1}{6u^2(s^*)}\sqrt{C_\Delta}\left(1+\frac{p(s^*)}{3u^2(s^*)}\right)
\nonumber\\
&=&-\mathbf{B}(s^*)_1.
\label{eq:B1}
\end{eqnarray}

Furthermore, by~\eqref{z},
\begin{eqnarray}
	\lefteqn{
	\lim_{s\downarrow s^*}
	\left(
\frac{\mathbf{\Psi}(s)_2-\mathbf{\Psi}(s^*)_2}{\sqrt{s-s^*}}
\right)
=
	\lim_{s\downarrow s^*}
	\left(
	\frac{z(s)-z(s^*)}{\sqrt{s-s^*}}
	\right)
}
\nonumber\\
&=&
	\lim_{s\downarrow s^*}
	\left(
\frac{1}{\sqrt{s-s^*}}
\left(\frac{\lambda x(s)}{2+2\lambda+2s-x(s)}
-\frac{\lambda x(s^*)}{2+2\lambda+2s^*-x(s^*)}\right)
\right)
\nonumber\\
&=&
\frac{2\lambda(1+\lambda+s^*)}{(2+2\lambda+2s^*-x(s^*))^2}
\mathbf{B}(s^*)_1
\nonumber\\
&=&-\mathbf{B}(s^*)_2,
\end{eqnarray}
which gives,
\begin{equation}
\lim_{s\downarrow s^*}
\left(
\frac{\mathbf{\Psi}(s)-\mathbf{\Psi}(s^*)}{\sqrt{s-s^*}}
\right)
=
-[\mathbf{B}(s^*)_1\ \ \ \mathbf{B}(s^*)_2]
=-{\bf B}(s^*),
\end{equation}
as expected~\eqref{eq:eqB}.
\\\bigskip

Now, assuming $\lambda=2.5$, we solve~\eqref{sstarex2} numerically,
\begin{eqnarray}
s^*&\approx& -1.1178,
\end{eqnarray}
and then evaluate $[x\ z]={\bf\Psi}(s^*)$ using~\eqref{eq:x} to get $x$ and then~\eqref{z} to get~$z$,
\begin{eqnarray}
{\bf\Psi}(s^*)&\approx&[1.7878\quad 1.5016],
\end{eqnarray}
and then ${\bf K}(s^*)$ and $(-{\bf D}(s^*))$ using~\eqref{eq:KsDs},
\begin{eqnarray}
{\bf K}(s^*)&\approx&[-2.0944],\\
-{\bf D}(s^*)&\approx&
\left[
\begin{array}{cc}
-0.5944 &   4.0016\\
2.0000 &  -0.8822
\end{array}
\right]
,
\end{eqnarray}
which have common eigenvalue $\gamma \approx -2.0944$. Also, we use~\eqref{Us*} to evaluate
\begin{eqnarray}
{\bf U}(s^*)&\approx& [3.5756\quad 3.0031].
\end{eqnarray}
Finally, we evaluate ${\bf B}(s^*)$ using~\eqref{eq:B1}, and ${\bf Y}(s^*)$ using~\eqref{eq1},
\begin{eqnarray}
{\bf B}(s^*)
&\approx&
[1.6416  \quad\quad\quad 2.2069],\\
{\bf B}(s^*){\bf Q}_{21}(s^*){\bf B}(s^*)
&\approx& [2.6948\quad\quad\quad 3.6228],\\
{\bf Y}(s^*)
&\approx&
[0.8808  \quad\quad -0.6197],
\end{eqnarray}
which gives
\begin{eqnarray}\label{eq:test}
{\bf K}(s^*){\bf B}(s^*)+{\bf B}(s^*){\bf D}(s^*)
&\approx &
10^{-14}\times [0.7550   \quad -0.0888],
\end{eqnarray}
which is approximately zero, as expected~\eqref{eq2}.

Finally,
\begin{eqnarray}
{\bf H}(s^*,y)&=&
\sum_{n=1}^{\infty}
\frac{y^n}{n!}
{\bf H}_{1,n}(s^*)
\nonumber\\
&=&
\sum_{n=1}^{\infty}
\frac{y^n}{n!}
\sum_{i=0}^{n-1}
\left({\bf K}(s^*) \right)^i
\times
\mathbf{B}(s^*){\bf Q}_{21}(s^*)
\times
\left({\bf K}(s^*) \right)^{n-1-i}
\nonumber\\
&=&
\sum_{n=1}^{\infty}
\frac{y^n}{n!}
\sum_{i=0}^{n-1}
\left( {\bf K}(s^*)\right)^{n-1}
\left(\mathbf{B}(s^*){\bf Q}_{21}(s^*)\right)
\nonumber\\
&=&
ye^{{\bf K}(s^*)y}
\mathbf{B}(s^*){\bf Q}_{21}(s^*)
\nonumber\\
&\approx&
1.6416 ye^{-2.0944\times y},
\end{eqnarray}
and
\begin{eqnarray}
{\bf H}(s^*)&=&
\int_{y=0}^{\infty}
ye^{{\bf K}(s^*)y}
{\bf B}(s^*){\bf Q}_{21}(s^*)
dy
\nonumber\\
&=&
({\bf K}(s^*))^{-2}
{\bf B}(s^*){\bf Q}_{21}(s^*)
,
\nonumber\\
\widetilde{\bf H}(s^*)&=&{\bf H}(s^*)
+
\left(
-({\bf K}(s^*))^{-1}\mathbf{B}(s^*)
+
{\bf H}(s^*){\bf\Psi}(s^*)
\right)
{\bf 1}
\nonumber\\
&\approx& 3.4428
,
\end{eqnarray}
and so
	\begin{eqnarray}
	\bmu(dy)^{(0)}_{11}&=& 
	diag({\widetilde{\bf H}(s^*)})^{-1}
	{ {\bf H}(s^*,y)}  
	dy
	\nonumber\\
	&\approx&
	0.2905\times
	1.6416 ye^{-2.0944\times y} 
	dy
	,\nonumber\\
	\bmu(dy)^{(0)}_{12}&=& 
	diag({\widetilde{\bf H}(s^*)})^{-1}
	\left(
	e^{{\bf K}(s^*)y}{\bf B}(s^*)
	+
	{\bf H}(s^*,y){\bf\Psi}(s^*)
	\right)
	dy
	\nonumber\\
	&\approx&
	0.2905\times
	\left(
	[2.6948\quad 3.6228]
	+1.6416 y
	[1.7878\quad 1.5016] 
	\right)	
	e^{-2.0944\times y}
	dy
	.
	\nonumber\\
	\end{eqnarray}

\end{example}

\bibliographystyle{abbrv}
\bibliography{all}
\end{document}